\documentclass{article}
\usepackage{latexsym}
\usepackage{amsmath}
\usepackage{amssymb}
\usepackage[dvips]{graphicx}
\usepackage{theorem}
\newcommand{\e}{{\rm e}}
\newcommand{\dive}{{\rm div}}
\newcommand{\erre}{{\mathbb R}}
\newcommand{\renne}{{\mathbb R}^n}
\newcommand{\rnpiu}{{\mathbb R}^n_+}
\newcommand{\enne}{{\mathbb N}}
\newcommand{\qed}{\hfill $\Box$}

\newcommand{\Om}{\Omega}
\newcommand{\sobdue}{W^{2,\infty}}

\newcommand{\dist}{{\rm dist}}

\newcommand{\esse}[2]{{\cal S}_{#2}^{#1}}
\newcommand{\io}{\hat{\text{\it \i}}}
\makeatletter
\@addtoreset{equation}{section}
\makeatother
\newtheorem{theorem}{Theorem}[section]
{\theorembodyfont{\rmfamily}\newtheorem{remark}[theorem]{Remark}}
\newtheorem{lemma}[theorem]{Lemma}
\newtheorem{prop}[theorem]{Proposition}
\newtheorem{definition}[theorem]{Definition}

\begin{document}

\null \vspace{1cm}
\begin{center}
{\LARGE Global calibrations for the non-homogeneous

\vspace{.5cm}

Mumford-Shah functional} \vspace{1cm}

 {\normalsize Massimiliano \sc Morini}

\vspace{.5cm}

{\normalsize S.I.S.S.A.}\\
{\normalsize via Beirut 2-4, 34014 Trieste, Italy}\\
{\normalsize e-mail:\tt\  morini@sissa.it}

\vspace{2cm}

\begin{minipage}[t]{11.5cm}
{\footnotesize \centerline{ {\bf Abstract}} \noindent Using a
calibration method we prove that, if $\Gamma\subset \Om$ is a
closed regular hypersurface and if the function $g$ is
discontinuous along $\Gamma$ and regular outside, then the
function $u_{\beta}$ which solves
$$
\begin{cases}
\Delta u_{\beta}=\beta(u_{\beta}-g)& \text{in
$\Om\setminus\Gamma$}\\
\partial_{\nu} u_{\beta}=0 & \text{on $\partial\Om\cup\Gamma$}
\end{cases}
$$
is in turn discontinuous along $\Gamma$ and it is the unique {\em
absolute minimizer} of the non-homogeneous Mumford-Shah functional
$$
\int_{\Om\setminus S_u}|\nabla u|^2\, dx +{\cal
H}^{n-1}(S_u)+\beta\int_{\Om\setminus S_u}(u-g)^2\, dx,
$$
over $SBV(\Om)$, for $\beta$ large enough. Applications of the
result to the study of the gradient flow by the method of minimizing movements
are shown.}\vspace{.5cm}\\


\vspace{.5cm} \noindent{\footnotesize {\bf Key words:}
free-discontinuity problems, calibration method, minimizing
movements}
\end{minipage}
\end{center}

\setcounter{page}{0} \thispagestyle{empty}

\vfill

\begin{center}
{\normalsize Ref. S.I.S.S.A. 41/2001/M (May 2001)}
\end{center}

\pagebreak

\clearpage

\null \vspace{1cm}
\begin{center}
{\LARGE Global calibrations for the non-homogeneous

\vspace{.5cm}

Mumford-Shah functional} \vspace{1cm}

 {\normalsize Massimiliano \sc Morini}

\vspace{1cm}

\begin{minipage}[t]{11.5cm}
{\footnotesize \centerline{ {\bf Abstract}} \noindent Using a
calibration method we prove that, if $\Gamma\subset \Omega$ is a
closed regular hypersurface and if the function $g$ is
discontinuous along $\Gamma$ and regular outside, then the
function $u_{\beta}$ which solves
$$
\begin{cases}
\Delta u_{\beta}=\beta(u_{\beta}-g)& \text{in
$\Omega\setminus\Gamma$}\\
\partial_{\nu} u_{\beta}=0 & \text{on $\partial\Omega\cup\Gamma$}
\end{cases}
$$
is in turn discontinuous along $\Gamma$ and it is the unique {\em
absolute minimizer} of the non-homogeneous Mumford-Shah functional
$$
\int_{\Omega\setminus S_u}|\nabla u|^2\, dx +{\cal
H}^{n-1}(S_u)+\beta\int_{\Omega\setminus S_u}(u-g)^2\, dx,
$$
over $SBV(\Omega)$, for $\beta$ large enough. Applications of the
result to the study of the gradient flow
by the method of minimizing movements are shown.}
\vspace{.5cm}\\
\end{minipage}
\end{center}
\section{Introduction}

The Mumford-Shah functional was introduced in \cite{Mum-Sha1}
within the context of a variational approach in Image
Segmentation. In the $SBV$ setting proposed by De Giorgi (see
\cite{DG-Amb}) it can be written as
$$
F(u)=\int_{\Om\setminus S_u}|\nabla u|^2\, dx +\alpha{\cal
H}^{n-1}(S_u)+\beta\int_{\Om\setminus S_u}(u-g)^2\, dx,
$$
where $g:\Om\to \erre$ is the given input function, $\alpha$ and
$\beta$ are positive parameters,
 ${\cal H}^{n-1}$ is the $(n-1)$-dimensional Hausdorff measure, $u$ is the unknown function
in  the space $SBV(\Om)$ of special functions of bounded variation
in $\Om$, $S_{u}$ is the set of essential discontinuity points of
$u$, while $\nabla u$ denotes its approximate gradient (see
\cite{Amb-Fus-Pal}).

In dimension two, the function $u$ which minimizes $F$ over
$SBV(\Om)$ (whose existence is stated in \cite{Amb}) can be
thought to represent a piecewise approximation of the input grey
level function $g$, while $S_u$ represents the set of relevant
contours in the image. One of the mathematical features of the
Mumford-Shah functional is a very strong lack of convexity, which
produces, for example, non-uniqueness of the solution and makes
the exhibition of explicit minimizers a very difficult task.
Concerning this last point, the calibration method recently
developed by Alberti, Bouchitt\'e, and Dal Maso in
\cite{Alb-Bou-DM} seems to be a powerful tool. For some
applications of this method see \cite{Alb-Bou-DM}, \cite{DM-M-M},
\cite{M-M}, or \cite{Circe}. Coming back to  $F$, throughout the paper we keep
the
 parameter $\alpha$ fixed (and, without loss of generality, equal
to 1) and we are interested in minimizers of the functional
\begin{equation}\label{fbetag}
F_{\beta,g}(u)=\int_{\Om\setminus S_u}|\nabla u|^2\, dx +{\cal
H}^{n-1}(S_u)+\beta\int_{\Om\setminus S_u}(u-g)^2\, dx,
\end{equation}
with  $g$ piecewise smooth function. It is intuitive that taking
$\beta$ large means penalizing a lot the $L^2$-distance between $g$
and the solution, which is therefore forced to be close to the
input function. More precisely it is easy to see that, if for
simplicity we take $g$ belonging to $SBV(\Om)$ such that
\begin{equation}\label{ipotesisug}
F_{\beta,g}(g)=\int_{\Om\setminus S_g}|\nabla g|^2\, dx +{\cal
H}^{n-1}(S_g)=C<+\infty,
\end{equation}
then, denoting by $u_{\beta}$  a minimum point of $F_{\beta,g}$,
we have
$$
\int_{\Om}(u_{\beta}-g)^2\, dx\leq
\frac{F_{\beta,g}(u_{\beta})}{\beta}\leq
\frac{F_{\beta,g}(g)}{\beta}=\frac{C}{\beta},
$$
that is $u_{\beta}\to g$ in $L^2(\Om)$  as $\beta\to+\infty$. This
suggests that, in agreement with our expectations, if $\beta$ is
large, then $u_{\beta}$ should be an accurate  reconstruction of
the original image $g$. Actually, T.J.Richardson in \cite{Rich}
has proved also the convergence of the discontinuity sets in
dimension two: more precisely, he has shown that if $g$ is a
function of class $C^{0,1}$ outside any neighbourhood of the
singular set $S_g$  satisfying \eqref{ipotesisug}, and if $S_g$
has no isolated points (i.e. for every $x\in S_g$ and for every
$\rho>0$, ${\cal H}^1(B_{\rho}(x)\cap S_g)>0$), then, as
$\beta\to+\infty$,
$$
S_{u_{\beta}}\to S_g\ \text{in the Hausdorff
metric}\qquad\text{and}\qquad {\cal H}^1(S_{u_{\beta}})\to{\cal
H}^1(S_{g}).
$$
In the main theorem of the paper (see Theorem \ref{calibration1}),
 using the calibration method mentioned above, we are able to
 prove that, under suitable assumptions on the regularity of $\Om$,
 $g$, and $S_g$, a much stronger result holds true:

 \noindent
 {\it Suppose that $\Gamma$ is a closed hypersurface of class
 $C^{2,\alpha}$ contained in the $n$-dimensional
 domain $\Om$ (satisfying in turn some regularity assumptions),
 and let $g$ a function belonging to
 $W^{1,\infty}(\Om\setminus\Gamma)$, with $S_g=\Gamma$ and
 $\inf_{x\in\Gamma}(g^+(x)-g^-(x))>0$ (where $g^+$ and $g^-$
 denote the upper and the lower traces of $g$ on $\Gamma$). Then
 there exists $\beta_0>0$ depending on $\Gamma$, on the
 $W^{1,\infty}$-norm of $g$, and on the size of the jump of $g$
 along $\Gamma$, such that, for $\beta\geq\beta_0$, $F_{\beta,g}$
 has a unique minimizer $u_{\beta}$ which satisfies
 $$S_{u_{\beta}}=\Gamma.$$}

Let  us give now a short insight into some technical aspects of
the proof; we start by recalling the theorem on which the
calibration method is based. We shall consider the collection ${\cal
F}(\Om\times\erre)$ of
 all bounded vector fields $\phi=(\phi^x,\phi^z):
 \Om\times\erre \to \erre^n { \times }
\erre$ with the following property: there
 exists a finite family $(U_i)_{i\in
I}$ of pairwise disjoint and Lipschitz open
 subsets of $\Om\times\erre$ whose
closures cover $\Om\times\erre$,
 and a family $(\phi_i)_{i\in I}$ of vector fields in $Lip(\overline{U_i},\erre^n { \times }
\erre)$ such that $\phi$ agrees at any point with one of the $\phi_i$.

An {\it absolute calibration \/} for $u\in SBV(\Om)$ in
$\Om\times\erre$ is a
 vector field $\phi \in {\cal
F}(\Om\times\erre)$ which satisfies the
 following properties:
\label{condcal}
\begin{description}
\item[(a)] $\dive \phi = 0$ in $U_i$, for every $i\in I$;
\item[(b)] $\nu_{\partial U_i}\cdot \phi^+=\nu_{\partial
U_i}\cdot\phi^-=\nu_{\partial U_i}\cdot\phi$ ${\cal H}^{n-1}$-a.e in
$\partial U_i$ for every $i\in I$, where $\nu_{\partial U_i}(x)$
denotes the (unit) normal vector at $x$ to $\partial U_i$,
while $\phi^+$ and $\phi^-$ denote the two traces  of $\phi$
on the two sides of $\partial U_i$;
\item[(c)]
$\displaystyle\frac{\left(\phi^x(x,z)\right)^2}{4} \leq
\phi^{z}(x,z)+\beta (z-g(x))^2$ for almost every
 $x\in\Om$
and every $z\in\erre$;

\item[(d)]  $\phi^{x}(x,u(x)) = 2 \nabla
u(x,y)$
 and $\phi^{z}(x,u(x)) = |\nabla u(x)|^{2}-\beta (g(x)-u(x))^2$
for almost every $x\in \Om \setminus S_u$;
\item[(e)] $ \displaystyle\int_{u^{-}(x)}^{u^{+}(x)} \phi^{x}(x,z)\, dz
= \nu_u (x)$ for ${\cal H}^{n-1}$-a.e. $x\in S_u$, where $\nu_u(x)$ denotes the
unit normal vector at $x$ to $S_u$, which points toward $u^+$;
\item[(f)] $\displaystyle \left| \int_{s}^{t} \phi^{x}(x,z)\, dz
\right| \leq 1$ for ${\cal H}^{n-1}$-a.e. $x\in \Om$ and
for every $s$, $t\in\erre$;
\item[(g)] $\phi^x(x,z)\cdot \nu (x)=0$ for ${\cal H}^n$-a.e.
$(x,z)\in\partial (\Om\times\erre)$, where $\nu (x)$ denotes  the unit normal
vector at $x$ to $\partial \Om$.
\end{description}
Note that conditions (a) and (b) imply that $\phi$ is divergence free
in the sense of distributions in $\Om\times\erre$.

The following theorem is proved in \cite{Alb-Bou-DM2}.

\begin{theorem}
If there exists an absolute calibration $\phi$ for $u$ in $\Om { \times }
\erre$, then $u$ is an  absolute minimizer of the
Mumford-Shah functional (\ref{fbetag}) over $SBV(\Om)$.
\end{theorem}

\begin{remark}\label{unicity}
If for a.e. $x\in\Om$ the inequality in (b) is strict for $z\neq u(x)$, then
$u$ is the {\em unique} absolute minimizer of \eqref{fbetag}. The proof can be
obtained arguing as in the last part of Paragraph 5.8 in \cite{Alb-Bou-DM2}.
\end{remark}

The main difficulty in constructing the calibration comes from the
fact that the candidate $u_{\beta}$, which is the solution of the
Euler equation
\begin{equation}\label{euler}
\begin{cases}
\Delta u_{\beta}=\beta(u_{\beta}-g) &\text{in
$\Om\setminus\Gamma$}\\
\partial_{\nu} u_{\beta}=0 &\text{on
$\partial(\Om\setminus\Gamma)$},
\end{cases}
\end{equation}
presents, in general, a non vanishing gradient and a nonempty
discontinuity set. We remark that the case of $g$ equal to
characteristic function of a regular set (i.e. with vanishing
gradient) and the case of $g$ regular in the whole $\Om$ (i.e.
with empty discontinuity set) have been already treated in
\cite{Alb-Bou-DM} and require  a  simpler construction. From the
point of view of calibrations, the interaction (actually the
clash) between the (non vanishing) gradient and the (nonempty)
discontinuity set is reflected in the fact that  we have to
guarantee simultaneously conditions (d) and (e), which push
in opposite directions. Indeed condition (d) says that $\varphi^x$
on the graph of $u$ is tangential  to $\Gamma$ while (e) implies
that $\varphi^x$ must be on the average orthogonal to $\Gamma$ for
$x\in\Gamma$ and $t$ between $u^-(x)$ and $u^+(x)$; so we have to
``rotate'' suitably $\varphi^x$,  preserving at the same time
condition (f). Another difficulty comes from the fact the we have
to estimate how quickly the gradient of $u_{\beta}$ changes
direction; indeed if near $\Gamma$ it becomes suddenly orthogonal
to $\Gamma$ and (e) holds true, it could happen  that condition
(f) is violated: this risk is overcome by carefully estimating the
 $L^{\infty}$-norm of the Hessian matrix
$\nabla^2u_{\beta}$ with respect to $\beta$. In
order to perform such an estimate we need to assume that $\Gamma$
is of class $C^{2,\alpha}$, for some $\alpha>0$. We underline
that, at least in dimension two, the regularity assumption is
close to optimal, since, by Bonnet Regularity Theorem (see
\cite{Bonnet}) (proved for $n=2$) in a neighbourhood of any
regular point the discontinuity set is of class $C^{1,1}$,
for every $g\in L^{\infty}(\Om)$.

The starting point for the construction of $\varphi$ is the
following remark: if we fiber a neighbourhood of the graph of
$u_{\beta}$ by the graphs of a family of function $(v_t)_{t\in\erre}$ all
satisfying \eqref{euler} and we let $\varphi$ to be the vector
field  equal to
$$
(2\nabla v_t,|\nabla v_t|^2-\beta(v_t-g)^2)
$$
on the graph of $v_t$, then $\varphi$ is divergence free.
Unfortunately this construction works only locally (i.e., in a
neighbourhood of $\Gamma$), but we will see how to modify it
suitably in order to get it working ``globally''.

As an application of our theorem, we give a proof of the following
fact: if $u_0$ is regular enough outside a smooth singular set
$S_{u_0}$, then the gradient flow $u(x,t)$ of $u_0$ (via
minimizing movements, see next section for the definition) for the
homogeneous functional
\begin{equation}\label{F0}
F_0(u)=\int_{\Om}|\nabla u|^2\, dx+{\cal H}^{n-1}(S_u),
\end{equation}
keeps, at least for small times,  the singular set of $u(\cdot,t)$
equal to $S_{u_0}$, while $u$ evolves in $\Om\setminus S_{u_0}$
according to the heat equation with Neumann boundary conditions on
$\partial(\Om\setminus S_{u_0})$. This result was proved in
dimension one by Gobbino (see \cite{Gob}), with a slightly
different definition of gradient flow.

The plan of the paper is the following. In Section 2 we recall some
definitions, fix some
notations, and collect some results which will be useful for the proof
of our theorems. In Section 3 we provide the estimates we mentioned above on
the norm of the solutions of \eqref{euler}. In Section 4 we give the proof
of the main result and, in dimension two, we extend it to the case of
$\Om$ with piecewise smooth boundary (say a curvilinear polygon) and of
$\Gamma$ touching the boundary (orthogonally). The final section is devoted
to the study of minimizing movements.

\section{Preliminary Results}

For fixed $R>0$, we introduce the following class of sets:
\begin{multline}\label{defU}
{\cal U}_{R}  =   \left\{  \right.   E \subset \erre^n, E \ \hbox{open} :
\forall p \in \partial E \  \exists \, p', p'':  \\
p \in \partial B(p', R) \cap \partial B(p'', R),
\left.  B(p', R) \subset E ,\, B(p'', R)\subset {\cal C}E\ \right\},
\end{multline}
and
\begin{equation}\label{defUomega}
{\cal U}_{R}(\Om)  =\left\{E\in {\cal U}_R:\,E\subset\Om,\,
\dist(E,\partial\Om)\geq R\right\}.
\end{equation} If $E$ belongs
to ${\cal U}_R$ and $p\in \partial E$, we denote the centers of
the interior and exterior balls associated with $p$ by $p'$ and
$p''$ respectively; moreover, we call $\esse{p}{E}$ the class of
all coordinate systems centred at $p$ such that the vector
$\frac{1}{2R} (p'' - p')$ coincides with the $n$-th vector of the
coordinate basis. The following  proposition is proved in
(\cite{mora-morini-99})
\begin{prop}\label{regular}
There exists a constant $\rho > 0$ (depending only on $R$), such that
for every $E\in {\cal U}_{R}(\Om)$ and
for every $p_0 \in \partial E$, if we call
$C$ the cylinder $\{ x\in \erre^{n-1} :|x| < \rho  \} { \times }
] { - } R,R[$ expressed with respect to a coordinate
system belonging to $\esse{p_0}{E}$,
then $\partial E \cap C$ is
the subgraph of a function $f$ belonging to
$\sobdue ( \{ x\in \erre^{n-1} :|x| < \rho  \}  )$.
Moreover, the $\sobdue$-norm of $f$ is bounded by a constant depending
only on $R$ (independent of $p_0$, of $E$ and of the choice of  the coordinate
system in
 $\esse{p_0}{E}$).
\end{prop}
\begin{remark}\label{rdefUomega}
Note that if $\Om$ is bounded and of class $C^2$ then there
exists $R>0$ such that $\Om\in {\cal U}_R$.
\end{remark}

For $E\subset\erre^n$, we define the signed distance function
$$
d_E(x)=\dist(x,E)-\dist(x,{\cal C}E).
$$
Now we are going to state some basic properties of that function; for a proof
see, for example, \cite{DZ}.
\begin{lemma}\label{-2}
\begin{description}
\item{i)} Let $x$ be a point of $\erre^n$. Then $d_E(x)$ is differentiable at
$x$ if and only if there exists a unique $y\in\partial E$
such that $|d_E(x)|=|x-y|$. In this case, we have
$$
\nabla d_E(x)=\frac{x-y}{d_E(x)}
$$
and we can define the projection on $\partial E$ $\pi_E(x):=y$.
\item{ii)} Let $\partial E$ be a hypersurface of class $C^k$,
$k\geq 2$. Then, for
every
 $x\in\partial E$, there exists a neighbourhood $V$ of $x$ such that
$d_E\in C^k(V)$ and $\pi_E\in C^{k-1}(V)$.
 \end{description}
\end{lemma}

\begin{lemma}\label{-1}
Let $E\subset\erre^n$ be an open set whose boundary is a hypersurface
of class $W^{2,\infty}$. Then for every $x\in\partial\Om$, there exists
a neighbourhood $V$ of $x$ where $\pi_E$ is well defined and such that
$d_{E}\in W^{2,\infty}(V(x))$. Moreover, denoting by
$\lambda_1\leq\dots\leq\lambda_n$ the eigenvalues of $\nabla^2d_E$
and by $k_1(y)\leq\dots\leq k_{n-1}(y)$ the principal curvatures of $\partial
E$ at $\pi(y)$, we have
$$
\lambda_i:=
\begin{cases}
0 &\text{if $i=1$}\\
\frac{k_{i-1}(y)}{1+ d_E(y)k_{i-1}(y)} &\text{if $i>1.$}
\end{cases}
$$
\end{lemma}

\begin{lemma}\label{0}
Let $E$ be an open set belonging to ${\cal U}_R$, for some $R>0$. Then the
projection $\pi_E$ is well defined and of class $W^{1,\infty}$ in the
$(R/2)$-neighbourhood of $\partial E$, and therefore $d_E$ is of class
$W^{2,\infty}$ in that neighbourhood. Moreover we have:
$$
\|d_E\|_{W^{2,\infty}}\leq
C\qquad\text{and}\qquad\|\pi_E\|_{W^{1,\infty}}\leq C ,
$$
where $C$ is a positive constant depending only on $R$.
\end{lemma}
{\sc Proof.} The fact that $\pi_E$ is well defined in the
$(R/2)$-neighbourhood of $\partial E$ (denoted by $(\partial
E)_{R/2}$ ) is an easy consequence of the definition of ${\cal
U}_R$: indeed let $x$ be a point of $(\partial E)_{R/2}\cap {\cal
C}E$ and let $p\in\partial E$ such that $d_E(x)=|x-p|$. We claim
that such a $p$ is unique. Indeed let $B(p'',R)\subset{\cal C}E$
be the exterior ball associated with $p$ (see the definition
\eqref{defU}); since the vector $p''-p$ is parallel to $x-p$
(indeed both vectors are normal to $\partial E$ at $p$), it is
clear that $\overline{B(x,d_E(x))}\setminus\{p\}\subset
B(p'',R)\subset{\cal C}E$ and so $p$ is the unique minimum point.

Concerning the smoothness, it is enough to prove that $d_E$ is of class
$W^{2,\infty}$, then we conclude by the equality
$$
\pi_E(x)=x-d_E(x)\nabla d_E(x).
$$
Exploiting the definition of ${\cal U}_R$ in a way similar to the
one we did above, we can easily see that, for every $\varepsilon\in
(0,R/2)$,
\begin{equation}\label{identity}
(E)_{\varepsilon}\in {\cal U}_{R-\varepsilon}\qquad\text{and}\qquad
d_{(E)_{\varepsilon}}=d_E-{\varepsilon},
\end{equation}
implying that $\partial((E)_{\varepsilon})$ is in turn of class
$W^{2,\infty}$. So if $x\in (\partial E)_{R/2}$, then $x\in
\partial((E)_{\varepsilon})$ for  $\varepsilon=d_E(x).$
By Lemma \ref{-1} there exists a neighbourhood $V$ of $x$ where
$d_{(E)_{\varepsilon}}$ is of class $W^{2,\infty}$ and
$\|d_{(E)_{\varepsilon}}\|_{W^{2,\infty}}\leq C$, with $C$
depending only on $R$. Recalling \eqref{identity}, we are
done.\qed

For the proof of the announced estimates on the norm
of the solutions of \eqref{euler},
we will use some technical results coming from sectorial
operators theory and from interpolation theory.

First let us recall what a sectorial operator is.

Let $X$ a complex Banach space and $A:D(A)\to X$ a closed linear operator with
not necessarily
 dense domain; call $\rho(A)$ the resolvent set of $A$ and for
$\lambda\in\rho(A)$ denote
 by $R(\lambda,A)$ the resolvent operator $(\lambda
I-A)^{-1}$ belonging to $L(X)$.
 \begin{definition}\label{sectorial}
$A$ is said to be sectorial (in $X$) if the following two conditions are satisfied:
\begin{description}
\item{i)} there exist $\omega\in\erre$ and $\theta\in(\frac{\pi}{2},\pi)$ such that
$$S_{\theta,\omega}:=\{\lambda\in {\mathbb C}:\, |{\rm arg}(\lambda-\omega)|\leq\theta\}\subset
\rho(A);$$
\item{ii)} there exists a positive constant $M$ such that, for every $\lambda\in S_{\theta,
\omega}$,  there holds
$$\|R(\lambda,A)\|_{L(x)}\leq \frac{M}{|\lambda-\omega|}.$$
\end{description}
\end{definition}
We recall that $D(A)$, endowed with the norm
$$\|x\|_{D(A)} = \|x\|_{X} + \|Ax\|_{X}$$
is a Banach space continuously embedded in $X$.

Let $\Om$ be either $\erre^n$ or $\rnpiu$ and let
$A:\Om\to\erre^{n\times n}$ be a  matrix with  coefficients belonging to
$W^{1,\infty}(\Om)$
and uniformly elliptic, i.e., satisfying
$$
A(x)\xi\cdot\xi\geq\lambda |\xi|^2\qquad \forall x\in \Om,\ \forall \xi\in\renne,
$$
where $\lambda_0$ is a suitable positive constant; set
$$
D({\cal A}_0):=\left\{u\in L^{\infty}(\Om):
\, u\in\bigcap_{p\geq 1}W_{loc}^{2,p}(\Om),\  \dive(A\nabla u)\in L^{\infty}(\Om)
 \hbox{ and }
 A\nabla u\cdot \nu = 0
\hbox{ on } \partial\Om\right\}
$$
$$D({\cal A}_1):=\left\{u\in D({\cal A}_0):\,  \dive(A\nabla u)\in W^{1,\infty}(\Om)
 \hbox{ and } A\nabla u\cdot \nu = 0
\hbox{ on } \partial\Om\right\},$$ where $\nu(x)$ denotes the outer
unit normal vector at $x$ to $\Om$,  and define the operators
\begin{equation}\label{A0}
\begin{array}{rcl}
{\cal A}_0:D({\cal A}_0) & \to & L^{\infty}(\Om) \\
u & \mapsto & f\dive(A\nabla u),
\end{array}
\end{equation}
and
\begin{equation}\label{A1}
\begin{array}{rcl}
{\cal A}_1:D({\cal A}_1) & \to & W^{1,\infty}(\Om) \\
u & \mapsto & f\dive(A\nabla u),
\end{array}
\end{equation}
where $f:\Om\to (0,+\infty)$ is a positive function of class
$W^{1,\infty}$ satisfying:
$$
f(x)\geq \lambda_1>0\qquad\forall x\in\Om.
$$
The following fact is proved in \cite{LU} (see Theorem 3.1.6, page 77,  Theorem
3.1.7, page 78, and 3.1.26, page 103).
\begin{theorem}\label{3.1.6}
The operators ${\cal A}_0$ and ${\cal A}_1$ are  sectorial in
$L^{\infty}(\Om)$ and $W^{1,\infty}(\Om)$ respectively. In particular there
exist two positive constants $\beta_0$ and $K$, depending on the constants
$\lambda_0$, $\lambda_1$, on
 $W^{1,\infty}$-norm of  $A$ and $f$,  such that the problem
\begin{equation}\label{NPV}
\begin{cases}
f\dive(A\nabla u) = \beta(u-g) & \text{in $\Om$}, \\
A\nabla u\cdot \nu = 0 & \text{in $\partial\Om$},
\end{cases}
\end{equation}
admits  a unique solution $u\in D({\cal A})$, for every
$\beta\geq\beta_0$ and for every $g\in L^{\infty}(\Om)$ . Moreover
$u$ satisfies
\begin{equation}\label{s0}
\|u\|_{\infty}+\beta^{-\frac{1}{2}}\|\nabla u\|_{\infty}\leq K\|g\|_{\infty};
\end{equation}
if $g$ belongs to $W^{1,\infty}(\Om)$ then the following estimate actually holds
\begin{equation}\label{s1}
\|u\|_{W^{1,\infty}} + \beta^{-\frac{1}{2}}\|f\dive(A\nabla u)\|_{\infty}
+ \sup_{x_0\in\Om}\beta^{\frac{n}{2p}}
\|\nabla^2 u\|_{L^p\left( B(x_0,\frac{1}{\sqrt\beta})\cap\Om\right)}
\leq K \|g\|_{W^{1,\infty}}.
\end{equation}
\end{theorem}

Given a sectorial operator $A:D(A)\to X$ there is a natural way to construct a family
of intermediate spaces between $D(A)$ and $X$, by setting for $\theta\in(0,1)$
$$
D(A,\theta,\infty) = \left\{x\in X:\ \sup_{t>2\omega\lor 1}\left(t^{\theta}
\|AR(t,A)x\|_{L(X)}\right)<+\infty\right\},
$$
where $\omega$ is the real number appearing in i) of Definition \ref{sectorial}.
Setting
\begin{equation}\label{sinterpol}
[x]_{D(A,\theta,\infty)} = \sup_{t>2\omega\lor 1}\left(t^{\theta}
\|AR(t,A)x\|_{L(X)}\right),
\end{equation}
one sees that $[x]_{D(A,\theta,\infty)}$ is a seminorm and
$D(A,\theta,\infty)$ endowed with the norm
\begin{equation}\label{ninterpol}
\|x\|_{D(A,\theta,\infty)} = \|x\|_X + [x]_{D(A,\theta,\infty)}
\end{equation}
is a Banach space. Moreover, for $0\leq\theta_1<\theta_2\leq 1$,
$$
Y\subseteq D(A,\theta_2, \infty)\subset D(A,\theta_1,\infty)\subseteq X,
$$
with continuous embeddings.
An important fact is stated in the following proposition
\begin{prop}\label{2.2.7} (see Proposition 2.2.7, page 50 of \cite{LU})
$$
\begin{array}{rcl}
A_{\theta}: D(A,\theta+1,\infty):=\{x\in D(A):\ Ax\in D(A,\theta,\infty)\}
& \to &  D(A,\theta,\infty)\\
x & \mapsto & Ax,
\end{array}
$$
is sectorial in $D(A,\theta,\infty)$; moreover
\begin{equation}\label{uholder}
\|R(\lambda, A_{\theta})\|_{L(D(A,\theta,\infty))}\leq
\|R(\lambda, A)\|_{L(X)}.
\end{equation}
\end{prop}
Next theorem gives a useful characterization of the intermediate spaces $D(A,\theta,\infty)$
in the case of elliptic operators.
\begin{theorem}\label{3.1.30}(see Theorem 3.1.30, page 108 of \cite{LU}) Let ${\cal A}_0$ be the
operator defined in \eqref{A0}. Then for every $\theta\in (0,\frac{1}{2})$,
$$
D({\cal A}_0,\theta,\infty)= C^{0,2\theta}(\overline{\Om}),
$$
with equivalence of the respective norms. In particular there exists two constants $C_1$ and
$C_2$ depending
only on the $W^{1,\infty}$-norm of $A$ and $f$ and on the constants
$\lambda_0$ and $\lambda_1$, such that
\begin{equation}\label{interpol}
C_1\|g\|_{D({\cal A}_0,\theta,\infty)}\leq
\|g\|_{C^{0,2\theta}(\overline{\Om})}\leq C_2\|g\|_{D({\cal
A}_0,\theta,\infty)}.
  \end{equation}
\end{theorem}

Let us recall now the definition of gradient flow
for the
homogeneous Mumford-Shah functional \eqref{F0} via
 minimizing movements  (see for instance
\cite{Cha-Dov} or \cite{Amb0}). Let $\Om$ be a bounded open subset of
$\erre^n$
 and consider an initial datum $u_0\in L^{\infty}(\Om)$. For fixed
$\delta>0$ (which is the time discretization parameter) we can
define the {\em $\delta$-approximate evolution}
$u_{\delta}(\cdot):[0,+\infty)\to SBV(\Om)$ as the affine
interpolation of the discrete function
$$
\begin{array}{rcl}
\delta\enne&\to& SBV(\Om)\\
\delta i&\mapsto& u_{\delta,i},
\end{array}
$$
where $u_{\delta,i}$ is
 inductively defined as follows:
 $u_{\delta,0}= u_0$ and $u_{\delta,i}$   is a solution of
 $$
\min_{v\in SBV(\Om)}\int_{\Om}|\nabla v|^2\, dx+{\cal
H}^{n-1}(S_v)+\frac{1}{\delta}\int_{\Om}|v-u_{\delta,i-1}|^2\, dx.
$$
  The existence of a solution of the problem above
is guaranteed by the Ambrosio theorem (see \cite{Amb}). We call {\em
minimizing movement} for $F_0$ with initial datum $u_0$, the set of
all functions $v:[0,+\infty)\to SBV(\Om)$ such that, for a
suitable subsequence $\delta_n\downarrow 0$, $u_{\delta_n}(t)\to
v(t)$ in $L^2(\Om)$, for every $t>0$.

\section{Technical Estimates}
\subsection{Estimates in smooth domains}\label{smooth}
Given a hypersurface
$\Gamma$ of class $C^{2,\alpha}$ we can define

\begin{equation}\label{normaholder}
\Lambda^{\alpha}({\Gamma}):=\sup_{x,y\in\Gamma}
\frac{\left|\nabla_{\tau}\nu(x)-\nabla_{\tau}\nu(y)\right|}
{|x-y|^{\alpha}},
\end{equation}
where $\nu$ is a smooth unit normal vector field to $\Gamma$
 and $\nabla_{\tau}$ denotes the tangential
 gradient along $\Gamma$.
\begin{lemma}\label{lholder}
Let $\Om$ be either $\renne$ or $\rnpiu$ and
 ${\cal A}_0$ be the operator  defined in \eqref{A0}. Then for every $\gamma\in
(0,\frac{1}{2})$ there exist two positive constants $K_0$ and
$\beta_0$, depending only on the constants of ellipticity
$\lambda_0$, $\lambda_1$, on $\gamma$, and on the
$W^{1,\infty}$-norm of the matrix $A$ and of the function $f$,
such that for every $\beta\geq \beta_0$ and for every $g\in
C^{0,1-\gamma}(\overline\Om)$ the solution $u$ of \eqref{NPV}
satisfies
\begin{equation}\label{sholder}
\beta^{\frac{1}{2}-\gamma}\|u-g\|_{C^{0,\gamma}(\overline\Om)}\leq
K_0 \|g\|_{C^{0,1-\gamma}(\overline\Om)}.
\end{equation}
\end{lemma}
{\sc Proof.} Recall that $u-g = {\cal A}_0R(\beta,{\cal A}_0)g$:
in order to obtain the thesis we have to estimate the quantity
$\displaystyle{\beta^{\frac{1}{2}-\gamma} \|{\cal
A}_0R(\beta,{\cal A}_0)g\|_{C^{0,\gamma}(\overline\Om)}}$. By
Theorems \ref{3.1.6} and \ref{3.1.30}, by \eqref{sinterpol} and
\eqref{ninterpol},
 there exist
$C_0>0$, $C_1>0$, and $\beta_0>0$, depending only on  $\lambda_0$,
$\lambda_1$, on $\gamma$,
 and on the $W^{1,\infty}$-norm of  $A$ and $f$,
such that
\begin{eqnarray}
\|{\cal A}_0R(\beta,{\cal A}_0)g\|_{C^{0,\gamma}(\overline\Om)}
&\leq& C_0
\|{\cal A}_0R(\beta,{\cal A}_0)g\|_{D({\cal A}_0,\frac{\gamma}{2},\infty)} \nonumber\\
& = & C_0\left(\|{\cal A}_0R(\beta,{\cal A}_0)g\|_{\infty}+
\sup_{t\geq2\beta_0\lor1}t^{\frac{\gamma}{2}} \|{\cal
A}_0R(t,{\cal A}_0){\cal A}_0R(\beta,{\cal
A}_0)g\|_{\infty}\right),\label{291}
\end{eqnarray}
and
\begin{equation}\label{293}
\sup_{2\beta_0\lor1\leq t} t^{\frac{1-\gamma}{2}}\|{\cal
A}_0R(t,{\cal A}_0)g\|_{\infty}\leq C_1
\|g\|_{C^{0,1-\gamma}(\overline\Om)}.
\end{equation}
We observe that \eqref{s0} implies the existence of two positive
constants $\beta_0$ and $C_2$, depending in turn on $\lambda_0$,
$\lambda_1$ and on the $W^{1,\infty}$-norm of $A$ and $f$, such
that
\begin{equation}\label{292}
\|\beta R(\beta,{\cal A}_0)\|_{L(L^{\infty}(\Om))}\leq C_2,
\end{equation}
for every $\beta\geq\beta_0$.
Using \eqref{292} and \eqref{293}, we can estimate
\begin{eqnarray}
\sup_{2\beta_0\lor1\leq\beta\leq
t}\beta^{\frac{1}{2}-\gamma}t^{\frac{\gamma}{2}} \|{\cal
A}_0R(t,{\cal A}_0){\cal A}_0R(\beta,{\cal A}_0)g\|_{\infty}
&=&\!\!\!\!\!\!\!\! \sup_{2\beta_0\lor1\leq\beta\leq
t}\beta^{\frac{1}{2}-\gamma}t^{\frac{\gamma}{2}}
\|{\cal A}_0R(\beta,{\cal A}_0){\cal A}_0R(t,{\cal A}_0)g\|_{\infty}\nonumber \\
& = &\!\!\!\!\!\!\!\!\sup_{2\beta_0\lor1\leq\beta\leq
t}\!\!\left(\frac{\beta}{t}\right)^ {\frac{1}{2}-\gamma}\!\!\!\!
t^{\frac{1-\gamma}{2}}
\|(\beta R(\beta,{\cal A}_0)-I){\cal A}_0R(t,{\cal A}_0)g\|_{\infty}\nonumber \\
&\leq& (C_2+1)\sup_{2\beta_0\lor1\leq t}
t^{\frac{1-\gamma}{2}}\|{\cal A}_0R(t,{\cal A}_0)g\|_{\infty}\nonumber \\
& \leq & (C_2+1)C_1\|g\|_{C^{0,1-\gamma}}, \label{294}
\end{eqnarray}
and analogously
\begin{equation}\label{295}
\sup_{2\beta_0\lor1\leq
t\leq\beta}\beta^{\frac{1}{2}-\gamma}t^{\frac{\gamma}{2}} \|{\cal
A}_0R(t,{\cal A}_0){\cal A}_0R(\beta,{\cal A}_0)g\|_{\infty} \leq
(C_2+1)C_1\|g\|_{C^{0,1-\gamma}}.
\end{equation}
Combining \eqref{294}, \eqref{295},  \eqref{291}, and using again \eqref{293}, we finally obtain
\begin{eqnarray*}
\sup_{\beta\geq2\beta_0\lor1}\beta^{\frac{1}{2}-\gamma} \|{\cal
A}_0R(\beta,{\cal A}_0)g\|_{C^{0,\gamma}(\overline\Om)} &\leq&
C_0\left(\sup_{\beta\geq2\beta_0\lor1}
\beta^{\frac{1}{2}-\gamma}\|{\cal A}_0R(\beta,{\cal A}_0)g\|_{\infty}\right.\\
& & \;\;\;\;\;\;\;+\left.
\sup_{\beta,t\geq2\beta_0\lor1}\beta^{\frac{1}{2}-\gamma}t^{\frac{\gamma}{2}}
\|{\cal A}_0R(t,{\cal A}_0){\cal A}_0R(\beta,{\cal A}_0)g\|_{\infty}\right) \\
&\leq& C_0(C_1+C_2+1)\|g\|_{C^{0,1-\gamma}}.
\end{eqnarray*}
\qed

The following theorem provides the preannounced estimate
on the Hessian $\nabla^2 u$ of the function $u$
which solves \eqref{euler}; we recall that $(\partial\Om')_{{R}}$
denotes the
${R}$-neighbourhood of $\partial\Om'$.
\begin{theorem}\label{casoreg}
Let $\Omega\subset\renne$ be a bounded domain of class $C^{1,1}$.
\begin{description}
\item{i)} For every $R>0$, we can
find two positive constants $\beta_0 = \beta_0(R)$ and $K = K(R)$
with the property that
 if $\Om'$ is a domain
belonging to ${\cal U}_R(\Om)$,  then
 for every $\beta\geq\beta_0$ and
 for every $g\in W^{1,\infty}\left(\Om\setminus\overline{\Om'}\right)$ the solution $u$ of
\begin{equation}\label{NP}
\begin{cases}
\Delta u = \beta(u-g) & \text{in $\Om\setminus\overline{\Om'}$,}\\
\partial_{\nu} u = 0 & \text{on
$\partial\left(\Om\setminus\overline{\Om'} \right)$,}
\end{cases}
\end{equation}
satisfies
\begin{equation}\label{scasoreg}
\|\nabla u\|_{\infty} + \beta^{-\frac{1}{2}}
\|\Delta u\|_{\infty} +
\beta^{\frac{n}{2p}-1}\sup_{x_0\in\Om\setminus\overline{\Om'}}
 \|\nabla^2
u\|_{L^p\left(B(x_0,\frac{1}{\sqrt\beta})\cap\Om\setminus\overline{\Om'}\right)}\leq K
 \|g\|_{W^{1,\infty}}.
\end{equation}
A similar  conclusion holds for the solution of
\begin{equation}\label{NP'}
\begin{cases}
\Delta u = \beta(u-g) & \text{in $\Om'$,}\\
\partial_{\nu} u = 0 & \text{on $\partial\Om'$}.
\end{cases}
\end{equation}
\item{ii)} For every $R>0$, for every $\overline{\Lambda}>0$,
and for every $\gamma\in(0,\alpha)$ (with $\alpha\in(0,1)$), there
exist two positive constants $\beta_0 =
\beta_0(R,\overline{\Lambda},\gamma)$ and $K =
K(R,\overline{\Lambda},\gamma)$ with the property that
 if $\Om'$ is a domain of class $C^{2,\alpha}$
 belonging to ${\cal U}_R(\Om)$,
and $\Lambda^{\alpha}
(\partial\Om')\leq\overline{\Lambda}$, then, for every $\beta\geq\beta_0$ and
 for every $g\in W^{1,\infty}\left(\Om\setminus\overline{\Om'}\right)$, the solution $u$ of
\eqref{NP}
satisfies
$$
\|\nabla^2u\|_{L^{\infty}\left((\partial\Om')_{{R}}\cap
(\Om\setminus\overline{\Om'})\right)}\leq K \beta^ {\frac{1}{2} +
\gamma}\|g\|_{W^{1,\infty}}.
$$
 A similar conclusion holds
for the solution of problem \eqref{NP'}.
\end{description}
\end{theorem}
{\sc Proof.} We will prove in details only {\it ii)}. Fix
$p\in\partial\Om'$. By Proposition \ref{regular} there exist two
positive constants $\eta$ and $M_1$, the former depending only on
$R$ while the latter also on $\Lambda^{\alpha}(\partial \Om)$ ,
such that the cylinder $C^{\eta}:=\{ x\in \erre^{n-1} :|x| < \eta
\} { \times } ] { - } R,R[$ (expressed with respect to a
coordinate system belonging to $\esse{p}{\Om'}$), intersected
with$\Om'$  is the subgraph of a function $f$ belonging to
$C^{2,\alpha} ( S )$ ($S:=C^{\eta}\cap\{x_n = 0\}$) and satisfying
\begin{equation}\label{sregular}
\|f\|_{C^{2,\alpha}}\leq M_1.
\end{equation}

Let $\theta\in C^{2,\alpha}_0(C^{\eta})$, $0\leq\theta\leq 1$ and
$\theta\equiv 1$ in $2^{-1}C^{\eta}$, such that
\begin{equation}\label{theta}
\partial_{\nu} \theta = 0\ \text{on $\partial\Om'\cap
C^{\eta}$}\qquad \text{and}\qquad \|\theta\|_{C^{2,\alpha}}\leq
M_2,
\end{equation}
where $M_2$ depends only on $R$.

Set $v = \theta u$ and note that $v$ solves
$$
\begin{cases}
\Delta v = \beta (v - h) & \text{in $\Om'\cap C^{\eta}$},\\
\partial_{\nu} v = 0 & \text{on $\partial(\Om'\cap
C^{\eta})$},
\end{cases}
$$
where $h:=\theta g + \beta^{-1}(\Delta\theta u + 2\nabla u\nabla\theta)$; finally, denoting
by $\psi$ the map
$$
\begin{array}{rcl}
C^{\eta} & \to & \psi(C^{\eta}) \\
(x_1,\dots,x_{n-1},x_n) & \mapsto &(x_1,\dots,x_{n-1},x_n-f(x_1,\dots,x_{n-1})),
\end{array}
$$
and setting $\tilde v := v\circ\psi^{-1}$ and $\tilde h := h\circ\psi^{-1}$, one sees
that (recall that $\tilde v$ and $\tilde h$ have compact support in $\psi(C^{\eta})$)
$$
\begin{cases}
\tilde f\dive(\tilde  A \nabla\tilde v) = \beta (\tilde v - \tilde h) & \text{in $\rnpiu$},\\
\tilde A\nabla\tilde v\cdot \nu = 0 & \text{on
$\partial(\rnpiu)$},
\end{cases}
$$
where $\tilde A$ and $\tilde f$ are $W^{1,\infty}$-extensions to
$\renne_+$ of the matrix-valued function $A :=
 \left[\frac{D\psi
(D\psi)^*}{|{\rm det}\psi|}\right]\circ\psi^{-1}$ and of the
function $f:=|{\rm det}\psi|\circ \psi^ {-1}$ respectively ,
satisfying
$$
\|\tilde A\|_{W^{1,\infty}(\renne_+)} =
\|A\|_{W^{1,\infty}(\psi(C^{\eta}))}, \qquad \|\tilde
f\|_{W^{1,\infty}(\renne_+)}=\|f\|_{W^{1,\infty}(\psi(C^{\eta}))}
$$
and
$$
\tilde A(x)\xi\cdot\xi\geq\frac{1}{2}|\xi|^2\ \forall x\in
\renne_+,\ \forall \xi\in\renne,\qquad \tilde f(x)\geq\frac{1}{2}\
\forall x\in\renne_+
$$
(since $A(0)= I$ and $f(0)=1$, by \eqref{sregular}, we can choose
$\eta$ depending only on $R$ such that the property above holds
true in $\psi(C^{\eta})$).

The solution $\tilde v$ can be suitably decomposed as $\tilde v =
\tilde v_1 + \tilde v_2+\tilde v_3$ in the following way: set $h_1
= \theta g$, $h_2 =  \beta^{-1} \nabla u\nabla\theta$,
$h_3:=\beta^{-1}\Delta\theta u$, and $\tilde h_i =
h_i\circ\psi^{-1}$ ($i=1,2,3$) and choose $\tilde v_i$ as the
solution of
$$
\begin{cases}
\dive(\tilde  A \nabla\tilde v_i) = \beta (\tilde v_i - \tilde h_i) & \text{in $\rnpiu$},\\
\tilde A\nabla\tilde v_i\cdot \nu = 0 & \text{on
$\partial(\rnpiu)$},
\end{cases}
$$
for $i=1,2,3$.

 Applying Lemma \ref{lholder} we have, for $i=1,2,3$,
\begin{equation}\label{301}
\beta^{\frac{1}{2}-\gamma}\|\tilde v_i-\tilde
h_i\|_{C^{0,\gamma}}\leq K_0 \|g\|_{C^{0,1-\gamma}},
\end{equation}
where $K_0$ is a constant depending only on $\gamma$ and on the
norm of $\tilde A$, therefore (by definition of $A$ and by
\eqref{sregular}) only on $\gamma$ and $R$.

{\it Estimate for $\tilde v_1$.}
>From \eqref{301}, \eqref{sregular}, \eqref{theta}, and the definition of $\tilde h_1$ we deduce
$$
\beta^{\frac{1}{2}-\gamma}\|\tilde v_1-\tilde
h_1\|_{C^{0,\gamma}}\leq K_0K_1
(\|g\|_{C^{0,1-\gamma}}+\beta^{-1}\|u\|_{C^{0,1-\gamma}}),
$$
where $K_1$ depends only on $R$, and therefore, since
 by  \eqref{uholder} and \eqref{interpol}, we have
 $$
\|u\|_{C^{0,1-\gamma}}\leq K_2 \|g\|_{C^{0,1-\gamma}},
 $$
we obtain
$$
\beta^{\frac{1}{2}-\gamma}\|\tilde v_1-\tilde
h_1\|_{C^{0,\gamma}}\leq K_0K_1K_2 \|g\|_{C^{0,1-\gamma}},
$$
where $K_2$ depends only on $R$.
Combining the above inequality with the well known Schauder estimate, we finally obtain
\begin{equation}\label{sv1}
\|\nabla^2\tilde v_1\|_{\infty}\leq K_3 \|\tilde f\dive(\tilde A \nabla
\tilde v_1)\|_{C^{0,\gamma}} =K_3
\beta^{\frac{1}{2}+\gamma}\beta^{\frac{1}{2}-\gamma} \|\tilde
v_1-\tilde h_1\|_{C^{0,\gamma}}\leq
K_3K_0K_1K_2\beta^{\frac{1}{2}+\gamma}
 \|g\|_{C^{0,1-\gamma}},
\end{equation}
where $K_3$ depends only on $C^{1,\gamma}$-norm of $A$ and $f$ and
therefore only on $R$ and $\overline{\Lambda}$.

{\it Estimate for $\tilde v_2$.}
Arguing exactly as in the previous point, we obtain
\begin{equation}\label{302}
\beta^{\frac{1}{2}-\gamma}\|\tilde v_2-\tilde
h_2\|_{C^{0,\gamma}}\leq K_0K_1 \beta^{-1}\|\nabla
u\|_{C^{0,1-\gamma}}.
\end{equation}
By the Sobolev Embedding Theorem and by estimate \eqref{scasoreg}
(with $p = \frac{n} {\gamma}$) we have, for $\beta\geq\beta_0$ and
for every $x\in\Om\setminus\overline{\Om'}$,
\begin{equation}\label{-}
[\nabla u]_{C^{0,1-\gamma}\left((\Om\setminus\overline{\Om'})\cap
B\left(x,\beta^{-\frac{1}{2}}\right)\right)}\leq
Q_0\|\nabla^2u\|_{L^{\frac{n}{\gamma}}(\Om\setminus\overline{\Om'}\cap
B(x,\beta^{-\frac{1}{2}}))} \leq
Q_0Q_1\beta^{1-\frac{\gamma}{2}}\|g\|_{W^{1,\infty}},
\end{equation}
and
\begin{equation}\label{--}
\|\nabla u\|_{\infty}\leq Q_1\|g\|_{W^{1,\infty}},
\end{equation}
where $Q_0$ is the constant of Sobolev Embedding and depends only
on $\gamma$ while $Q_1$ depends only on $R$. If $|x-y|\geq
\beta^{-\frac{1}{2}}$, then, by \eqref{--}, we infer
\begin{equation}\label{---}
\frac{|\nabla u(x)-\nabla u(y)|}{|x-y|^{1-\gamma}}\leq
\beta^{\frac{1-\gamma}{2}} 2\|\nabla u\|_{\infty}\leq
2Q_1\beta^{\frac{1-\gamma}{2}}\|g\|_{W^{1,\infty}}.
\end{equation}
Combining \eqref{-}, \eqref{--}, and \eqref{---}, we get
$$
\|\nabla u\|_{C^{0,1-\gamma}}\leq
Q_1(Q_0+1)\beta^{1-\frac{\gamma}{2}} \|g\|_{W^{1,\infty}},
$$
and substitution in \eqref{302}, together with  Shauder Estimate,
yields
\begin{equation}\label{sv2}
\|\nabla^2\tilde v_2\|_{\infty}\leq K_3 \|\tilde f\dive(\tilde A
\nabla \tilde v_2)\|_{C^{0,\gamma}} = K_3
\beta^{\frac{1}{2}+\gamma}\beta^{\frac{1}{2}-\gamma} \|\tilde
v_2-\tilde h_2\|_{C^{0,\gamma}}\leq K_3K_0K_1Q_1(Q_0+1)
\beta^{\frac{1+\gamma}{2}}
 \|g\|_{W^{1,\infty}}.
\end{equation}

{\it Estimate for $\tilde v_3$.} First we note that, by
\eqref{uholder} and \eqref{interpol},
$$
\|\tilde v_3\|_{C^{0,\gamma}}\leq K_4 \|\tilde
h_3\|_{C^{0,\gamma}},
$$
with $K_4$ depending only on $R$; so we can estimate
\begin{eqnarray*}
\|\tilde v_3-\tilde h_3\|_{C^{0,\gamma}}&\leq& \|\tilde
v_3\|_{C^{0,\gamma}}+\|\tilde h_3\|_{C^{0,\gamma}}\\
&\leq& (K_4+1)\|\tilde h_3\|_{C^{0,\gamma}}\leq \beta^{-1}(K_4+1)M
\|u\|_{C^{0,\gamma}}\leq\beta^{-1}(K_4+1)K_4M
\|g\|_{W^{1,\infty}}.
\end{eqnarray*}
By Shauder Estimate we finally obtain,
\begin{equation}\label{sv3}
\|\nabla^2\tilde v_3\|_{\infty}\leq K_3 (K_4+1)K_4\|g\|_{W^{1,\infty}}.
\end{equation}

By  \eqref{sregular} and again \eqref{scasoreg} we  have
\begin{eqnarray*}
\|\nabla^2u\|_{L^{\infty}(2^{-1}C^{\eta})}&\leq& C
\left(\|\nabla^2\tilde
v\|_{L^{\infty}(\rnpiu)}+ \|\tilde
v\|_{W^{1,\infty}(\rnpiu)}\right)\\
&\leq &CC' \left(\|\nabla^2\tilde
v_1\|_{L^{\infty}(\rnpiu)}+\|\nabla^2\tilde v_2\|_{L^{\infty}(\rnpiu)}
+\|\nabla^2\tilde v_3\|_{L^{\infty}(\rnpiu)}+
\|g\|_{W^{1,\infty}(\rnpiu)}\right),
\end{eqnarray*}
where $C$ and $C'$ depend only on $R$. Using \eqref{sv1},
\eqref{sv2}, and \eqref{sv3}, we finally deduce for
$\beta\geq\beta_0\lor 1$
$$
\|\nabla^2u\|_{L^{\infty}(2^{-1}C^{\eta})}\leq
CC'C''\left(\beta^{\frac{1}{2}+\gamma}
\|g\|_{W^{1,\infty}(\rnpiu)}+\|g\|_{W^{1,\infty}(\rnpiu)}\right)\leq2CC'C''
\beta^{\frac{1}{2}+\gamma}\|g\|_{W^{1,\infty}(\rnpiu)},
$$
where $C''$ depends only on $\gamma$, $R$, and
$\overline{\Lambda}$. Repeating all the above argument for every
$p\in\partial\Om'$ we get {\it ii)}.

The proof of statement {\it i)} can be done in a similar way: by
localizing, straightening  the boundary, and using Theorem
\ref{3.1.6}. \qed

\subsection{Estimates in domains with angles}

In the following $\Om\subset\erre^2$ will denote  a curvilinear
polygon which means that $\partial\Omega$ is given by the union of
a finite number of simple connected curves $\tau_1,\dots,\tau_k$
of class $C^3$ (up to their endpoints) meeting at corners with
different angles $\alpha_j\in (0,\pi)$ ($j=1,\dots,k$). Finally we
will denote by ${\cal S}$ the set of the vertices, i.e. the set of
the singular points of $\partial\Om$.

\begin{prop}\label{prop1}
Let $\Om$ be as above. Then there exists $\beta_0>0$ and $K>0$
such that for every $\beta>\beta_0$ and for every $g\in
L^{\infty}(\Om)$, the solution $u$ of
\begin{equation}\label{ubeta}
\begin{cases}
\Delta u=\beta(u-g)&\text{in $\Om$}\\
\partial_{\nu}u=0 &\text{on $\partial\Om$},
\end{cases}
\end{equation}
satisfies
\begin{equation}\label{sgri}
\|u\|_{\infty}+\beta^{-\frac{1}{2}}\|\nabla u\|_{\infty}\leq K
\|g\|_{\infty}.
\end{equation}
\end{prop}
{\sc Proof.} The  estimate is proved in \cite{GRI1} for the
corresponding Dirichlet problem in a polygon, but one easily sees
that the same proof actually works also in our case: indeed the
change of boundary conditions does not affect the argument, and
the main tool, which is a Calderon-Zygmund type inequality, proved
in \cite{GRI2}, is actually available also for curvilinear
polygon, as shown, for example, in \cite{Mum-Shah}.\qed

The following proposition is proved in \cite{Mum-Shah}.

\begin{prop}\label{prop2}
Let $\Om$ be as above. Then there exists $K>0$ such that for every $\beta>0$
$\beta>0$ and for every $g\in W^{1,\infty}(\Om)$, the function $u$
solution
of \eqref{ubeta}, satisfies:
\begin{equation}\label{sMum-Shah}
\beta^{\frac{1}{2}}\|u-g\|_{\infty}\leq K \|\nabla g\|_{\infty}.
\end{equation}
\end{prop}

\begin{prop}\label{prop3}
Let $\Om$ be as above. Then there exists a positive constant $K$
such that for every $\beta\geq 1$ and for every $g\in
W^{1,\infty}(\Om)$,
 the solution $u$  of \eqref{ubeta} satisfies:
\begin{equation}\label{gradirreg}
\|\nabla u\|_{\infty}\leq K
\|g\|_{W^{1,\infty}}\beta^{\frac{1}{4}}.
\end{equation}
\end{prop}
{\sc Proof.} Fix $\beta\geq 1$; by Proposition \ref{prop1} there
exists $\lambda_0>0$ independent of $\beta$ such that, setting
$g_{\lambda}=\frac{\Delta u-\lambda u}{\lambda}$, for
$\lambda\geq\lambda_0$ we have
\begin{equation}\label{interpol2}
\|\nabla u\|_{\infty}\leq K\sqrt{\lambda}\|g_{\lambda}\|_{\infty}\leq
K\sqrt{\lambda}\left(\frac{\|\Delta u\|_{\infty}}{\lambda}+
\|u\|_{\infty}\right)=K\left(\frac{\|\Delta
u\|_{\infty}}{\sqrt{\lambda}}+
\sqrt{\lambda}\|u\|_{\infty}\right).
\end{equation}
Now set $\lambda_{\rm min}:=\frac{\|\Delta
u\|_{\infty}}{\|u\|_{\infty}}$ and suppose that $\|\Delta
u\|_{\infty}\geq \lambda_0\|g\|_{\infty}$. It follows that
$\lambda_{\rm min}\geq \lambda_0$ (recall that $\|u\|_{\infty}\leq
\|g\|_{\infty}$):  therefore, taking $\lambda=\lambda_{\rm min}$
in \eqref{interpol2}, we obtain
$$
\|\nabla u\|_{\infty}\leq 2K\|\Delta u\|_{\infty}^{\frac{1}{2}}
\|u\|_{\infty}^{\frac{1}{2}}
$$
and therefore, by Proposition \ref{prop2},
$$
\|\nabla u\|_{\infty}\leq 2K\|g\|_{\infty}^{\frac{1}{2}}
\left(K'\beta^{\frac{1}{2}}\|\nabla
g\|_{\infty}\right)^{\frac{1}{2}}\leq K''
\|g\|_{W^{1,\infty}}\beta^{\frac{1}{4}},
$$
where $K''$ is independent of $\beta$.

If $\|\Delta u\|_{\infty}<\lambda_0\|g\|_{\infty}$, then we simply
use the Calderon-Zygmund type estimate proved in \cite{Mum-Shah}
(it is crucial here the hypothesis that all the angles are less
than $\pi$) to get the existence of a constant $C>0$, depending
only on $\Om$, such that
$$
\|u\|_{W^{2,p}}\leq C\|g\|_{\infty}\leq
C\|g\|_{\infty}\beta^{\frac{1}{4}}. $$ We conclude by applying the
Sobolev Embedding Theorem. \qed

\begin{prop}\label{BOOO}
Let $\Om$ and ${\cal S}$ be as above and let $\Gamma$ be a simple
connected curve in $\Om$ joining two points $x_1$ and $x_2$ belonging to
$\partial\Om\setminus{\cal S}$. Suppose in addition that $\Gamma$
is of class $C^3$ up to $x_1$ and $x_2$ (actually it would be
enough to take $\Gamma$ of class $C^3$ in two neighbourhoods $U_1$
and $U_2$ of $x_1$ and $x_2$ respectively, and of class
$C^{2,\alpha}$, for some $\alpha>0$, outside those neighbourhood).
Let us call $\Om_1$ and $\Om_2$ the two connected components of
$\Om\setminus\Gamma$. Finally set $\overline d:=\dist(x_1,{\cal
S})\land \dist(x_2,{\cal S})$. Then for every $\delta<\overline d$
and $\gamma\in\left(0,\frac{1}{2}\right)$, there exist two
positive constants $\beta_0$ and $K$ depending on $\delta$,
$\gamma$,  and $\Gamma$, such that, for every $\beta\geq\beta_0$
and for every $g\in W^{1,\infty}(\Om_i)$ ($i=1,2$), the solution
$u_{i}$ of
\begin{equation}\label{ubetai}
\begin{cases}
\Delta u_{i}=\beta(u_{i}-g)&\text{in $\Om_i$}\\
\partial_{\nu}u_i=0 &\text{on $\partial\Om_i$},
\end{cases}
\end{equation}
satisfies
\begin{equation}\label{scasoirreg}
\|\nabla u_{i}\|_{L^{\infty}\left((\Gamma)_{\delta}\cap
\Om_i\right)}+\beta^{-\left(\frac{1}{2}+\gamma\right)}
\|\nabla^2u_i\|_{L^{\infty}
\left((\Gamma)_{\delta}\cap\Om_i\right)}\leq
K\|g\|_{W^{1,\infty}}.
\end{equation}
\end{prop}
{\sc Proof.} The estimate can be performed by a localization
procedure as for Theorem \ref{casoreg} and in fact we have only to
look at what happens in a neighbourhood of $x_1$ and $x_2$. We
will look only at $x_1$ considered as a point of $\partial\Om_1$,
the other cases being analogous.

 First of all, as in \cite{Mum-Shah}, we  can
find  a neighbourhood $U=B(x_1,r)\cap\Om_1$ of $x_1$, for a
suitable $r\leq\delta$, and a diffeomorphism which transforms $U$
into a right angle, more precisely
 we can construct a one-to-one map
$\Phi=(\Phi_1,\Phi_2):U\cap\Om_1\to \Phi(U\cap\Om_1)$ of class
$C^{1,1}$ such that $\nabla\Phi(0,0)=I$ and
$\Phi(U)=\{w=(w_1,w_2)\in\erre^2:\, w_1>0,\, w_2>0 \}\cap V$,
where $V$ is a neighbourhood of the origin; we can endow $\Phi$
with the further property that if $v$ is a function defined in $U$
with normal derivative vanishing on $\partial\Om\cap
\overline{U}$, then $v\circ\Phi^{-1}$ has normal derivative
vanishing on $\Phi(\partial\Om\cap\overline{U})$ and vice-versa.
It follows, in particular, that $\Phi_2(x)$ has the following
properties:
\begin{itemize}
\item $\Phi_2(x)=0$ for every $x\in\Gamma\cap U$;
\item $\partial_{\nu}\Phi_2=0$ on $\partial\Om\cap \overline{U}$.
\end{itemize}
It is easy to see that we can choose a positive convex function $f$ such that
$$
f(0)=0,\ f'(0)=0,\ \text{and}\ \Delta(f\circ \Phi_2)\geq 0\text{
on $U':= B(x_1,r')\cap\Om_1$,}
$$
with $r'\leq r$, if needed. Thus we see that $f\circ \Phi_2$ is a
subsolution of
$$
\begin{cases}
\Delta u=0 &\text{in $U'$}\\
u=0 &\text{on $\Gamma\cap \overline{U'}$}\\
\partial_{\nu}u=0 &\text{on $\partial \Om\cap \overline{U'}$}\\
u=f\circ \Phi_2 &\text{on $\partial
U'\setminus(\partial\Om\cup\Gamma)$}
\end{cases}
$$
and therefore $f\circ \Phi_2\leq u$ in $U'$. By Theorem 5.1.3.1 of
\cite{GRI2} (actually it is stated only for polygons, but it can
be extended to curvilinear polygons, by the continuity method
used, for example, in \cite{Mum-Shah}) and the Sobolev Embedding
Theorem, $u$ is in $C^2(\overline{U''})$, where
$U''=B(x_1,r'')\cap\Om$, with $r''<r'$. Therefore, since $\nabla
(f\circ \Phi_2)(x_1)\neq 0$, and so $\nabla u(x_1)\neq 0$, we can
say that the map $\Psi:=(v,u)$, where $v$ is the harmonic
anticonjugate of $u$, is conformal in a neighbourhood
$U''':=B(x_1,r''')\cap\Om_1$, with $r'''\leq r''$, it belongs to
$C^2(\overline {U'''})$ and $\Psi(U''')=
\{w=(w_1,w_2)\in\erre^2:\, w_1>0,\, w_2>0 \}\cap V$, where $V$ is
a neighbourhood of the origin. Now take a cut-off function
$\theta$ of class $C^3$ such that $\theta\equiv 1$ on
$B(x_1,r'''/2)\cap\Om_1$, $\theta(x)=0$ for $|x|\geq (2/3)r'''$,
and $\partial_{\nu}\theta=0$ on $\partial\Om\cup\Gamma\cap
\overline{U'''}$; note that $v_{1}:=(\theta u_{1})\circ\Psi^{-1}$
solves
$$
\begin{cases}
A(w)\Delta v_{1}=\beta(v_{1}-h)&\text{in $\Psi(U''')$ }\\
\partial_{\nu}v_{1}=0 &\text{on
$\{w_1=0\}\cup\{w_2=0\}\cap\overline{\Psi(U''')}$},
 \end{cases}
$$
where $h:= [\theta g+\beta^{-1}(\Delta\theta u+2\nabla
u\nabla\theta)]\circ\Psi^{-1}$ and $A:=|\nabla u|^2\circ\Psi^{-1}$.

Moreover we have that $\partial_{\nu} A=0$ on $\{w_1=0\}\cap
\overline{\Psi(U''')}$, indeed, in view of the conformality of
$\Psi$, this is equivalent to say that $\partial_{\nu}|\nabla
u|^2=0$ on $\partial\Om\cap \overline{U'''}$, which is true by the
following computation
$$
\partial_{\nu}|\nabla u|^2=\partial_{\nu}(\partial_{\tau} u)^2=
2\partial_{\tau}u\partial^2_{\nu\tau}u=0,
$$
where we used the fact that $u\in C^2(\overline {U'''})$ and
$\partial_{\nu}u\equiv 0$ on $\partial\Om\cap U'''$. As a
consequence, the function
$$
\tilde A:=
\begin{cases}
A(w_1,w_2) &\text{if $w_1>0$ and $(w_1,w_2)\in\Psi(U''')$}\\
A(-w_1,w_2)&\text{if $w_1<0$ and $(-w_1,w_2)\in\Psi(U''')$}
\end{cases}
$$
turns out to be of class $C^1$ up to the boundary; in particular
it can be extended to a function, still denoted by $\tilde A$,
belonging to $C^{1}(\overline{\erre^2_+})\cap
W^{1,\infty}(\erre^2_+)$. Now it is easy to check that, denoting
by $\tilde v_{1}$ and $\tilde h$ the extensions by reflection of
$v_{1}$ and $h$ respectively,
$$
\begin{cases}
\tilde A(w)\Delta \tilde v_{1}=\beta(\tilde v_{1}-\tilde
h)&\text{in
$\erre^2_+$ }\\
\partial_{\nu}\tilde v_{1}=0 &\text{on $\{w_2=0\}$};
 \end{cases}
$$
at this point we are in a position to apply the regularity
theorems stated in Subsection \ref{smooth}, obtaining the desired
estimate for $\tilde v_{1}$. To complete the proof we can now
proceed exactly as we did for Theorem \ref{casoreg}.\qed

\section{The calibration}

\subsection{The regular case}
Let $\Om\subset\erre^n$ be  a bounded open subset of class
$C^{1,1}$ and let $\Om_1\subset \Om$ be an open set belonging to
${\cal U}_R(\Om)$ (see \eqref{defUomega}). We set
$\Om_2:=\Om\setminus\overline{\Om_1}$, $\Gamma:=\partial \Om_1$,
and, for every $x\in\Gamma$, we denote  the unit
outer normal to $\partial \Om_1$ at $x$ by $\nu(x)$.

\begin{lemma}\label{derivnorm}
There exist two positive constants $c$ and $\beta_0$, depending
only on $R$, such that, for every $\beta\geq \beta_0$, we can find
two functions $z_{1,\beta}:\Om_1\to\erre$ and
$z_{2,\beta}:\Om_2\to\erre$ of class $W^{2,\infty}$ with the
following properties:
\begin{description}
\item{i)}$\frac{1}{2}\leq z_{i,\beta}\leq 1$ in $\Om_i$, for $i=1,2$ and
$z_{2,\beta}\equiv \frac{1}{2}$ in a neighbourhood of $\partial
\Om$;
\item{ii)}$\Delta z_{i,\beta}\leq c\beta z_{i,\beta}$ in $\Om_i$, for
$=1,2$;
\item{iii)}$z_{1,\beta}(x)=z_{2,\beta}(x)=1$ and $\partial_{\nu}
z_{1,\beta}(x)=-\partial_{\nu}
 z_{2,\beta}(x)\geq\sqrt\beta$ for every
$x\in\Gamma$;

\item{iv)} $\|\nabla z_{i,\beta}\|_{\infty}\leq c\sqrt{\beta}$
and  $\|\nabla^2
z_{i,\beta}\|_{\infty}\leq c{\beta}$.
\end{description}
\end{lemma}
{\sc Proof.} Let us denote the signed distance function
from $\Om_1$ by $d$ and let $\pi$ the projection on $\Gamma$ which, by
Lemma \ref{0}, is well defined in $(\Gamma)_{\frac{R}{2}}$; we
begin by constructing $z_{2,\beta}$. Let
$w_{\beta}:[0,+\infty)\to(0,+\infty)$ be the solution of the
following problem
$$
\begin{cases} w_{\beta}''=16\beta w_{\beta},\\
w_{\beta}(0)=1/{2}, \\
w_{\beta}'({R}/{2})= 0,
\end{cases}
$$
which can be explicitly computed and it is given by
\begin{equation}\label{explicit}
w_{\beta}(t)=\frac{1}{2}\frac{\e ^{-4\sqrt{\beta}\frac{R}{2}}}{\e
^{4\sqrt{\beta}\frac{R}{2}}+\e ^{-4\sqrt{\beta}\frac{R}{2}}}\e
^{4\sqrt{\beta}t}+\frac{1}{2}\frac{\e
^{4\sqrt{\beta}\frac{R}{2}}}{\e ^{4\sqrt{\beta}\frac{R}{2}}+\e
^{-4\sqrt{\beta}\frac{R}{2}}}\e ^{-4\sqrt{\beta}t},
\end{equation}
 and let $\theta:[0,+\infty)\to[0,1]$ be a $C^{\infty}$ function such that
\begin{equation}\label{242}
\theta\equiv1\quad\text{in $[0,R/4]$}\qquad \theta\equiv 0\quad
\text{in $[R/2,+\infty)$} \qquad\text{and}\qquad
\|\theta\|_{C^2}\leq c_0,
\end{equation}
with $c_0$ depending only on $R$. We are now ready to define
$z_{2,\beta}:\Om_2\to\erre$ as
$$
z_{2,\beta}(x):=
\begin{cases}
\theta(d(x))((w_{\beta}(d(x))+1/2)+(1-\theta(d(x))){1}/{2} &
\text{if
$0<d(x)\leq R/2$},\\
{1}/{2} & \text{otherwise in $\Om_2$}.
 \end{cases}
$$
First of all note that, as it is  a convex combination of two
functions with range contained in $[1/2,1]$, $z_{2,\beta}$ itself
has range in $[1/2,1]$. Using the expression in \eqref{explicit} it is easy to
see that there exist
 $\beta_0>1$ and $c_1>1$ depending only on $R$ such that
\begin{equation}\label{241}
w_{\beta}'(0)\leq -\sqrt\beta\qquad |w'_{\beta}|\leq
c_1\sqrt\beta\quad\text{in
$[0,R/2]$}\qquad\text{and}\qquad|w''_{\beta}|\leq
c_1\beta\quad\text{in $[0,R/2]$},
\end{equation} for every $\beta\geq\beta_0$. From the first
inequality we obtain immediately {\it iii)} for $z_{2,\beta}$.
 Moreover, by \eqref{242} and  \eqref{241}, we can estimate
\begin{eqnarray*}
|\nabla z_{2,\beta}|&=& |\theta(d)w_{\beta}'(d)\nabla
d+\theta'\nabla d\, w_{\beta}(d)|\\
&\leq& |w_{\beta}'|+|\theta'|\leq c_1\sqrt\beta+c_0\leq
c\sqrt\beta,
\end{eqnarray*}
with $c$ depending only on $R$. Finally, using again
\eqref{242},\eqref{241}, and Lemma \ref{0}, we have
\begin{eqnarray*}
|\nabla^2 z_{2,\beta}|&\leq&|w_{\beta}'||\nabla^2
d|+|w_{\beta}'||\theta'|+|w_{\beta}''|+
 |\theta''|+|\theta'||\nabla^2 d|
+|\theta'||w_{\beta}'| \\
&\leq&c_1c_2\sqrt\beta+c_0c_1\sqrt\beta+c_1\beta+c_0+c_0c_2
+c_0c_1\sqrt\beta,
\end{eqnarray*}
where all the constants depend only on $R$ so that we can state
the existence of $c>0$, still depending only on $R$, such that
$$
|\nabla^2 z_{2,\beta}|\leq c \beta\qquad\forall\beta\geq\beta_0.
$$
 To conclude, we define
$z_{1,\beta}:\Om_1\to\erre$ as follows:
$$\qquad\qquad\quad
z_{1,\beta}(x):=
\begin{cases}
\theta(-d(x))((w_{\beta}(-d(x))+1/2)+(1-\theta(-d(x))){1}/{2} &
\text{if
$0>d(x)\geq -R/2$},\\
{1}/{2} & \text{otherwise in $\Om_1$}.\qquad\qquad\text{\qed}
 \end{cases}
$$

\begin{theorem}\label{calibration1}
Let $\Om\subset\erre^n$ be  a bounded open set of class
$C^{1,1}$ and let $\Om_1\subset \Om$ be an open set of class
$C^{2,\alpha}$ for some $\alpha\in(0,1)$ and compactly
contained in $\Om$. Let $R>0$
such that $\Om_1\in{\cal U}_R(\Om)$ (see \eqref{defUomega} and
 Remark \ref{rdefUomega})and set
$\Gamma:=\partial \Om_1$. Then for every function $g$ belonging
$W^{1,\infty}(\Om\setminus\Gamma)$, discontinuous along $\Gamma$
(i.e., $S_g=\Gamma$) and such that $g^+(x)-g^-(x)>S>0$ for every
$x\in\Gamma$, there exists $\beta_0>0$ depending on $R$, $S$,
 $\Lambda^{\alpha}(\Gamma)$ (see \eqref{normaholder}), and
$\|g\|_{W^{1,\infty}}$, such that for $\beta\geq\beta_0$ the
solution $u_{\beta}$ of
\begin{equation}\label{ubeta'}
\begin{cases}
\Delta u_{\beta}=\beta(u_{\beta}-g) &\text{in
$\Om\setminus\Gamma$,}\\
\partial_{\nu} u_{\beta}=0 &\text{on
$\partial\Om\cup\Gamma$},
\end{cases}
\end{equation}
is discontinuous along $\Gamma$ ($S_{u_{\beta}}=\Gamma$) and it is
the unique absolute minimizer of $F_{\beta,g}$ over $SBV(\Om)$.
\end{theorem}
{\sc Proof.}
 In the sequel we will denote the
signed distance from $\Om_1$ by $d$ and
 the projection on $\Gamma$ by $\pi$: by
Lemma \ref{0}, the two functions are well defined
 in
$(\Gamma)_{R/2}$. Moreover, in that neighbourhood, $d$ and $\pi$
are at least of class $W^{2,\infty}$ and
 $W^{1,\infty}$ respectively.

 As announced in the Introduction, the proof will be performed by
 constructing a calibration $\phi$;
adopting  the notation  introduced there, the vector field $\phi$
will be written as
$$
\phi(x,z) = (\phi^x(x,z),\phi^z(x,z)),
$$
where $\phi^x(x,z)$ is a $n$-dimensional ``horizontal'' component, while $\phi^z$ is
the (one dimensional) ``vertical'' component.
\begin{itemize}
\item {\it Preparation.}
\end{itemize}
Without loss of generality we can suppose that
$g^+$ coincides with the trace on $\Gamma$
of $g$ from $\Om_1$,
while $g^-$ is trace from $\Om_2$.
First of all let us choose $\beta'$, depending only on $R$, $S$,
and $\|g\|_{W^{1,\infty}}$  and $G$ depending  on $R$, such that,
for $\beta\geq\beta'$,
\begin{equation}\label{vicine}
\|u_{\beta}-g\|_{{\infty}}\leq \frac{S}{16}
\quad\text{and}\quad\sqrt\beta\|u_{\beta}-g\|_{{\infty}}\leq
G\|g\|_{W^{1,\infty}}\quad i = 1,2:
\end{equation}
this is possible by virtue of Theorem \ref{casoreg}.

As a second step, it is convenient to extend the
restriction of $u_{\beta}$
to $\Om_i$
($i=1,2$) to a $C^{1,1}$ function $u_{i,\beta}$ defined in
the whole $\Om$, in such a way that
\begin{equation}\label{utilde}
u_{i,\beta}(x) = u_{\beta}(x) \ \text{for
$x\in{\Om}_i$,}\quad \|
u_{i,\beta}\|_{W^{2,\infty}} \leq c
\|u_{\beta}\|_{W^{2,\infty}}, \quad\text{and}\quad
u_{1,\beta}-u_{2,\beta} \geq \frac{3}{4}S\quad \text{for
every $x\in \Om$},
\end{equation}
where $c$  is a positive constant depending only on $R$: this
operation can be performed in many ways, for example, to construct
$u_{2,\beta}$ we can extend the resctriction of $u_{\beta}$
to $\Om_2$ in a
neighbourhood of $\Gamma$ by a standard localization procedure and
then we can make a convex combination through a cut-off function
with $u_{\beta}-(3/4)S$ (recall that by definition of $S$ and by
\eqref{vicine}, we have $u^+_{\beta}- u^-_{\beta}>(3/4)S$ on
$\Gamma$); it is clear that all can be done in such a way that the
constant $c$ depends only on the ``$C^{1,1}$-norm'' of $\Gamma$
and therefore only on $R$. We require also that
$$
\partial_{\nu}  u_{1,\beta}=0\qquad\text{on $\partial\Om$.}
$$

 By \eqref{scasoreg} and \eqref{utilde}, we can
state the existence of two positive constants $K$ and $\beta''$
depending only on $R$  such that
\begin{equation}\label{grad}
\|\nabla u_{i,\beta}\|_{\infty}\leq
K\|g\|_{W^{1,\infty}}\qquad i=1,2,
\end{equation}
for every $\beta\geq \beta''$.

Let $\beta'''>0$ satisfying
\begin{equation}\label{beta''}
\frac{1}{6}\sqrt{\beta'''}=\max\left\{4(K\|g\|_{W^{1,\infty}})^2,
{64}/{S^2},\beta',\beta'',\beta_0\right\}+1,
\end{equation}
where $\beta_0$ is the constant appearing in  Lemma
\ref{derivnorm}. Let $z_{1,\beta'''}$ and $z_{2,\beta'''}$ be the
two functions constructed in Lemma \ref{derivnorm} with
$\lambda=\beta'''$ and define $v_{1}$, $v_{2}$ as follows
$$
v_{1}(x)=
\begin{cases}
z_{1,\beta'''}(x) & \text{if $x\in\overline{\Om_1}$}\\
2-z_{2,\beta'''}(x) &\text{if $x\in\Om_2$}
\end{cases}
$$
and
$$
v_{2}(x)=
\begin{cases}
z_{2,\beta'''}(x) & \text{if $x\in\overline{\Om_2}$}\\
2-z_{1,\beta'''}(x) &\text{if $x\in\Om_1$}.
\end{cases}
$$
From the properties of $z_{i,\beta}$ ($i=1,2$), as stated in Lemma
\ref{derivnorm}, it follows immediately that $v_{i}\in
W^{2,\infty}(\Om)$ and
\begin{equation}\label{sagmon}
\|\nabla v_{i}\|_{\infty}\leq K_1 \sqrt{\beta'''},\qquad\qquad
\|\nabla^2 v_i\|_{\infty}\leq K_1\beta'''\quad i=1,2
\end{equation}
where $K_1$ is a positive constant depending only on $R$. Note
that $\nabla v_1(x)=-\nabla v_2(x)$ for every $x\in\Om$. We
remark also that, for $x\in\Gamma$, by construction,
\begin{equation}\label{normale}
\frac{\nabla v_1(x)}{|\nabla v_1(x)|}=
-\frac{\nabla v_2(x)}{|\nabla v_2(x)|}=\nu(x),
\end{equation}
where $\nu(x)$ denotes the unit normal vector at $x$ to $\Gamma$
(outer with respect to $\Om_1$). We set
\begin{equation}\label{tildeh}
\tilde h_(x) = \frac{1}{\sqrt{2}}
|\nabla v_1|^{-\frac{1}{2}}=\frac{1}{\sqrt{2}}|\nabla v_2|^{-\frac{1}{2}}
\end{equation}
for every $x\in\Gamma$. Moreover, using \eqref{sagmon}
and {\it iii)} of Lemma \ref{derivnorm},
we can find a positive constant $D\leq R/2$, depending only on $R$,
$S$,
and $\|g\|_{W^{1,\infty}}$, such that
\begin{equation}\label{Dtilde}
|\nabla v_i(x)|\geq\frac{1}{2},\quad
\tilde h^2(\pi(x))\frac{|\nabla v_i(x)|}{v_i(x)}<1-
\frac{25}{32}\frac{1}{\sqrt 3}
\qquad\text{if $|d(x)|\leq D$,
$i=1,2$.}
\end{equation}
 Applying {\it iii)} of  Lemma \ref{derivnorm}, we get
\begin{equation}\label{pre-h}
|\nabla v_{i}(x)| ^{\frac{1}{2}}\geq
\sqrt[4]{\beta'''}\geq \max\left\{{8}/{S},1\right\}\qquad i = 1,2,
\end{equation}
 where the last inequality follows directly from \eqref{beta''}.

Moreover, combining Lemma \ref{derivnorm}, \eqref{grad}, and
\eqref{beta''}, we deduce
\begin{equation}\label{v12}
4|\nabla
u_{i,\beta}(x)|^2-\frac{1}{6}|\nabla v_{i}(x)|\leq
4(K\|g\|_{W^{1,\infty}})^2-\frac{1}{6}\sqrt{\beta'''}\leq -1
\end{equation}
and analogously
\begin{equation}\label{v12bis}
\frac{1}{\sqrt2}|\nabla v_i(x)|^{-\frac{1}{2}}\|\nabla
u_{i,\beta}\|_{\infty}<\frac{1}{4\sqrt3}\qquad i=1,2,
 \end{equation}
for every $x\in\Gamma$ and for every $\beta\geq\beta'''$.

Let $\varepsilon\in (0,1)$ be  such that
\begin{equation}\label{epsilon}
6\varepsilon\|\nabla u_{i,\beta}\|_{\infty}+4\varepsilon^2
\|\nabla v_i\|_{\infty}\leq\frac{1}{4}\quad \text{for $i = 1,2$
and $\beta\geq \beta'''$};
\end{equation}
by \eqref{grad} and \eqref{sagmon} (and the definition of
$\beta'''$) we see
 that $\varepsilon$ can be chosen depending only $R$,
$S$ and $\|g\|_{W^{1,\infty}}$.  By \eqref{tildeh},  it follows,
for every $x\in\Gamma$,
$$
4(\tilde h)^2\|\nabla v_i\|_{\infty}\geq 4(\tilde
h)^2(-1)^{i+1}\partial_{\nu} v_i= 4\cdot
\frac{1}{2}>\frac{1}{4},
$$
therefore, by \eqref{epsilon},
\begin{equation}\label{51}
\varepsilon< \tilde h(x)\quad \forall x\in\Gamma.
\end{equation}

Let $\gamma$ be a fixed constant belonging to
$(0,\frac{1}{2}\land\alpha)$: by applying {\it ii)} of Theorem
\ref{casoreg}, we can find two positive constants $\beta^{\text{\i
v}}$ and $K_2$  depending only on $R$ and
$\Lambda^{\alpha}(\Gamma)$ (and $\gamma$) such that
\begin{equation}\label{hessian}
\|\nabla^2
u_{\beta}\|_{L^{\infty}\left((\Gamma)_{\frac{R}{2}}\right)}\leq
K_2 \beta^{\frac{1}{2}+\gamma}\|g\|_{W^{1,\infty}},
 \end{equation}
for every $\beta\geq\beta^{\text{\i v}}$.

We can define, for $\beta>0$,
$$
 h_{\beta}(x) =
\begin{cases}
\left(\tilde
h(\pi(x))-\beta^{\frac{1}{2}+\gamma_1}|d(x)|\right)\lor
\varepsilon & \text{if $|d(x)|\leq  D$}\\
 \varepsilon & \text{if $|d(x)|> D$},
\end{cases}
$$
where $\gamma_1$ is a fixed constant belonging to
$(\gamma,\frac{1}{2})$. It is easy to see that there exists
$\beta^{\rm v}>0$ depending  on $D$ (and therefore only on $R$, $S$, and
$\|g\|_{W^{1,\infty}}$) such that
 $h_{\beta}$ is continuous
(in fact Lipschitz) for $\beta>\beta^{\rm v}$.

Using \eqref{pre-h}, \eqref{tildeh}, \eqref{sagmon}, and Lemma
\ref{0}, we have
\begin{equation}\label{gradh}
\|\nabla h_{\beta}\|_{\infty}\leq C'\left(\frac{1}{\sqrt
2}\left(\frac{S}{8}+1\right)^4
 \|\nabla^2 v_i\|_{\infty}
\|\nabla \pi\|_{\infty} +\beta^{\frac{1}{2}+\gamma_1}
\|\nabla
 d\|_{\infty} \right)\leq K_3 \beta^{\frac{1}{2}+\gamma_1},
\end{equation}
where $K_3$ is a positive constant depending on $R$, $S$, and
$\|g\|_{W^{1,\infty}}$.

Finally we set
\begin{equation}\label{beta1}
\beta_1 = \max\{\beta'',\beta''',\beta^{\text{\i v}},\beta^{\rm v
},1\}
\end{equation}
and
\begin{equation}\label{mui}
\mu_i(x)=\frac{\Delta v_i(x)}{v_i(x)};
\end{equation}
notice that by \eqref{sagmon} we get
\begin{equation}\label{mu}
\mu_i(x)\leq \frac{K_1\beta'''}{v_i(x)}\leq
2K_1\beta_1\qquad\text{for every $x\in\Om$ }.
\end{equation}
\begin{itemize}
\item{\it Definition of the calibration.}
\end{itemize}
From now on we will assume $\beta\geq \beta_1$. Let us consider
the following sets
\begin{equation}\label{Ai}
A_i := \{(x,z)\in\Om\times\erre:\
u_{i,\beta}(x)-h_{\beta}(x)\leq z\leq  u_{i,\beta}(x) +
h_{\beta}(x)\}, \quad i = 1,2.
\end{equation}
Since, by \eqref{utilde},$  u_{1,\beta}(x)-
u_{2,\beta}(x)\geq \frac{3}{4}S$ everywhere, noting that
$h_{\beta}\leq S/8$ everywhere (by \eqref{pre-h} and
\eqref{tildeh}), we  see
 that
$$
\dist(A_1, A_2)\geq \frac{S}{2}\qquad\text{for $\beta\geq
\beta_1$}.
$$
The crucial  point is in  constructing  the vector field  around
the graph of $u_{\beta}$, i.e. in $A_i\cap(\Om_i\times\erre)$:
here we have to provide a divergence free vector field satisfying
condition (d) of Section 1 and such that
\begin{eqnarray*}
\phi^x(x,z)\cdot \nu_{u_{\beta}}\geq 0 && \text{for $x\in\Gamma$ and
$u_{2,\beta}<z<u_{1,\beta}$,}\\
\phi^x(x,z)\cdot \nu_{u_{\beta}}\leq 0 && \text{for $x\in\Gamma$,
$z<u_{2,\beta}$ or $z>u_{1,\beta}$.}
\end{eqnarray*}
These properties are crucial in order to obtain (e) and (f) simultaneously.

The remaining work is  a matter of finding a suitable extension
which preserves all the properties of calibrations.

We start by giving the global definition of the horizontal component $\phi^x$:
\begin{equation}\label{horizontal}
\phi^x(x,z): =
\begin{cases}
2\nabla u_{i,\beta}-2\frac{
u_{i,\beta}-z}{v_i}\nabla v_i -\frac{16}{h_{\beta}}\left(
(-1)^i(z-u_{i,\beta})-\frac{h_{\beta}}{2}\right)^+\nabla
{u}_{i,\beta}& \text{if $(x,z)\in A_i$, $i=1,2$,}\\
0 & \text{otherwise in $\Om\times\erre$.}
\end{cases}
\end{equation}
Concerning $\phi^z$, we begin by defining it in
$A_i\cap\left(\overline\Om_i\times\erre\right)$:
\begin{equation}\label{verticalparz}
\phi_i^z(x,z):=\left|\nabla
u_{\beta}-\frac{u_{\beta}-z}{v^i}\nabla v^i\right|^2-
 \beta(z-g)^2  +(\beta-\mu_i)(u_{\beta}-z)^2
+ \Psi_{i}(x,z) \quad\forall (x,z)\in
A_i\cap\left(\overline\Om_i\times\erre\right),
\end{equation}
 where
$\mu_i$ is the function defined in \eqref{mui} and
$$
\Psi_{i}(x,z):=\int_{ u_{i,\beta}}^z
\dive_x\left[\frac{16}{h_{\beta}}\left(
(-1)^i(t-u_{i,\beta})-\frac{h_{\beta}}{2}\right)^+\nabla u_{i,\beta}
 \right]\,
 dt.
$$
Let us clarify that in the formulas above $(\cdot)^+$ stands for
$(\cdot)\lor 0$.

For $x\in\Om_i$ and $-h_{\beta}<(-1)^i(z-u_{\beta})<\frac{h_{\beta}}{2}$,
the field $\phi$  reduces to
\begin{equation}\label{reduces}
 \left(2\nabla u_{\beta}-2\frac{u_{\beta}-z}{v_i}\nabla v_i,\,
 \left|\nabla u_{\beta}-
\frac{u_{\beta}-z}{v_i}\nabla v_i\right|^2-
\beta(z-g)^2+(\beta-\mu_i)(u_{\beta}-z)^2\right)
\end{equation}
and so, by some easy computation and using the definition of
$u_{\beta}$ and $\mu_i$, we have
\begin{eqnarray*}
\dive\phi(x,z) &=& 2\left(\Delta u_{\beta} - \frac{u_{\beta}-z}{v_i}\Delta v_i\right)
-2\beta(z-g)-2(\beta-\mu_i)(u_{\beta}-z) \\
&=&2\beta(u_{\beta}-g) -
2\mu_i(u_{\beta}-z)-2\beta(z-g)-2(\beta-\mu_i)(u_{\beta}-z) =0.
\end{eqnarray*}
For $x\in\Om_i$ and $\frac{h_{\beta}}{2}<(-1)^i(z-u_{\beta})<h_{\beta}$,
$\phi$ is the sum of the field in \eqref{reduces} and
$$
\left(\text{$-\frac{16}{h_{\beta}}\left(
(-1)^i(z-u_{\beta})-\frac{h_{\beta}}{2}\right)^+\nabla {u_{\beta}},\,
\Psi_i(x,z)$}\right),
$$
which is clearly  divergence free by the definition of $\Psi_i$.
Eventually we have,
\begin{equation}\label{aparz}
\dive\phi=0\quad\text{ in }(\Om_i\times\erre)\cap A_i.
\end{equation}
It is time now to extend the definition of $\phi^z$. Before writing
the explicit expression, we remark that conditions (a) and (b)
of Section 1 imply that such extension is essentially unique.
More precisely, if $(U_j)_{j=1,\dots,10}$ is the family of all connected
components of $(\Om\times\erre)\setminus(\partial A_1\cup\partial A_2\cup
(\gamma\times\erre))$, it easy to see that $\phi^z$ is uniquely determined on
$(\Om\setminus\Gamma)\times\erre=\cup_{j=1}^{10}U_j$ by \eqref{horizontal},
\eqref{verticalparz}, and the two following necessary conditions:
\begin{itemize}
\item $\partial_z\phi^z=-\dive_x\phi^x$ in $U_j$ for $j=1,\dots,10$ (which
ensures condition (a) of Section 1),
\item $\phi^+\cdot\nu_{\partial
U_j}=\phi^-\cdot\nu_{\partial U_j}$ on $\partial U_j$ for every
$j=1,\dots,10$, where $\phi^+$ and $\phi^-$ denote the traces of $\phi$
on the two sides of $\partial U_j$.
\end{itemize}
The only freedom is in the choice of $\phi^z$ on $\partial U_j$
according to the condition
$$
\phi\cdot\nu_{\partial U_j}=\phi^+\cdot\nu_{\partial
U_j}=\phi^-\cdot\nu_{\partial U_j}.
$$
We are now ready to give the complete the definition of $\phi^z$;
for $(x,z)\in(\Om_1\times\erre)\setminus A_1$ we define $\phi^z(x,z)$ as
follows: \begin{equation}\label{vertical1}
\begin{cases}
\phi^x(x, u_{\beta}+h_{\beta})\cdot (-\nabla u_{\beta}-\nabla h_{\beta})+
\phi^z(x,u_{\beta}+h_{\beta})&\text{if $z>u_{\beta}+h_{\beta}$,}\\
\\
\phi^x(x, u_{\beta}-h_{\beta})\cdot (-\nabla u_{\beta}+\nabla h_{\beta})+
\phi^z(x,u_{\beta}-h_{\beta})&\text{if
$u_{\beta}-h_{\beta}>z\geq u_{2,\beta}+h_{\beta}$,}\\
\\
\chi_1(x,z)+\phi^z(x,u_{2,\beta}+h_{\beta})
+\phi^x(x, u_{2,\beta}+h_{\beta})\cdot (\nabla u_{2,\beta}+\nabla h_{\beta})&
\text{if $u_{2,\beta}+h_{\beta}>z\geq u_{2,\beta}-h_{\beta}$,}\\
\\
\phi^x(x, u_{2,\beta}-h_{\beta})\cdot (-\nabla u_{2,\beta}+\nabla h_{\beta})+
\phi^z(x,u_{2,\beta}-h_{\beta})&\text{if
$u_{2,\beta}-h_{\beta}>z$,}
 \end{cases}
\end{equation}
where
$$
\chi_1(x,z)=\int_z^{u_{2,\beta}+h_{\beta}}\dive_x
\phi^x(x,t)\, dt.
$$
We remark that in first and in the second line we used
the definition of $\phi^z$ already given in \eqref{verticalparz}, in the third
line we used the definition of $\phi^z(x, u_{2,\beta}+h_{\beta})$ given in the
second one, and finally in the last line we exploited the definition $\phi^z(x,
u_{2,\beta}-h_{\beta})$ given in the previous one.

Analogously, for $(x,z)\in(\Om_2\times\erre)\setminus A_2$ we define
$\phi^z(x,z)$ as follows:
\begin{equation}\label{vertical2}
\begin{cases}
\phi^x(x, u_{\beta}-h_{\beta})\cdot (-\nabla u_{\beta}+\nabla h_{\beta})+
\phi^z(x,u_{\beta}-h_{\beta})&\text{if $z<u_{\beta}-h_{\beta}$,}\\
\\
\phi^x(x, u_{\beta}+h_{\beta})\cdot (-\nabla u_{\beta}-\nabla h_{\beta})+
\phi^z(x,u_{\beta}+h_{\beta})&\text{if
$u_{\beta}+h_{\beta}<z\leq u_{1,\beta}-h_{\beta}$,}\\
\\
\chi_2(x,z)+\phi^z(x,u_{1,\beta}-h_{\beta})
+\phi^x(x, u_{1,\beta}-h_{\beta})\cdot (\nabla u_{1,\beta}-\nabla h_{\beta})&
\text{if $u_{1,\beta}-h_{\beta}<z\leq u_{1,\beta}+h_{\beta}$,}\\
\\
\phi^x(x, u_{1,\beta}+h_{\beta})\cdot (-\nabla u_{1,\beta}-\nabla h_{\beta})+
\phi^z(x,u_{1,\beta}+h_{\beta})&\text{if
$u_{1,\beta}+h_{\beta}<z$,}
 \end{cases}
\end{equation}
where
$$
\chi_2(x,z)=\int_z^{u_{1,\beta}-h_{\beta}}\dive_x
\phi^x(x,t)\, dt.
$$
Finally we set
$$
\phi^z(x,z)=0\qquad\text{on $(\Gamma\cap\erre)\setminus(A_1\cup A_2)$;}
$$
this concludes the definition of $\phi$ which, by construction (and recalling
\eqref{aparz}), satisfies conditions (a) and (b) of Section 1.

\begin{itemize}
\item $\phi^z+\beta(z-g)^2>|\phi^x|^2/{4}$ {\it for almost
every} $(x,z)\in\Om\times\erre$ {\it with} $z\neq u(x)$.
\end{itemize}
We first prove the condition above in $A_i\cap (\Om_i\times\erre)$, and then
in the remaining.
For $x\in\Om_i$ and $-h_{\beta}\leq (-1)^i(z-u_{\beta})\leq \frac{h_{\beta}}{2}$, by
\eqref{reduces}, we have that
$$
\phi^z+\beta(z-g)^2=\frac{|\phi^x|}{4}^2+(\beta-\mu_i)(u_{\beta}-z)^2>
\frac{|\phi^x|}{4}^2,
$$
so condition (c) of Secton 1 is trivially satisfied, with
strict inequality.

For $x\in\Om_i$ and
$\frac{h_{\beta}}{2}<(-1)^i(z-u_{\beta})\leq h_{\beta}$, using the
definition of $\phi$ we see that (c) is equivalent to
\begin{equation}\label{equivalent}
(1):=(\beta-\mu_i)(u_{\beta}-z)^2+\Psi_i(x,z)>
\frac{(16)^2}{4(h_{\beta})^2}|\nabla
u_{\beta}|^2([\cdots])^2-\frac{16}{h_{\beta}}[\cdots]\left(\nabla
u_{\beta}-\frac{u_{\beta}-z}{v_i}\nabla v_i\right) \nabla
u_{\beta}:=(2),
\end{equation}
where we wrote $[\cdots]$ instead of
$\left[((-1)^i(z-u_{\beta})-\frac{h_{\beta}}{2})^+\right]$; by \eqref{51},
\eqref{vicine}, \eqref{grad},
and \eqref{gradh}, we have
\begin{eqnarray*}
\Psi_{i}(x,z)
 &\geq& \int_{u_{\beta}}^{z}\left(-\frac{16}{h_{\beta}}[\cdots]|\Delta
u_{\beta}|
   - |\nabla u_{\beta}|\cdot
\left|\nabla\left(\frac{16}{h_{\beta}}[\cdots]\right)\right|\right)\, dt\\
&\geq& -\frac{16}{\varepsilon}S^2|\beta( u_{\beta}-g)|- S|\nabla
u_{\beta}| \left(\frac{16}{\varepsilon}(|\nabla u_{\beta}|+|\nabla
h_{\beta}|)+
\frac{16S}{\varepsilon^2}|\nabla h_{\beta}|\right)\\
&\geq& -\frac{16}{\varepsilon}S^2\sqrt\beta
G\|g\|_{W^{1,\infty}}-SK\|g\|_{W^{1,\infty}}
\left(\frac{16}{\varepsilon}(K\|g\|_{W^{1,\infty}}+K_3\beta^{\frac{1}{2}+\gamma_1})+
\frac{16S}{\varepsilon^2}K_3\beta^{\frac{1}{2}+\gamma_1}\right),
\end{eqnarray*}
therefore, recalling that the all the constants appearing in the
last expression depend only on $R$, $S$,  and
$\|g\|_{W^{1,\infty}}$, there exists a positive constant $C$
depending on the same quantities such that
\begin{equation}\label{spsi}
\Psi_{i}(x,z)\geq -C \beta^{\frac{1}{2}+\gamma_1};
\end{equation}
recalling that $|u_{\beta}-z|\geq
\frac{h_{\beta}}{2}\geq\frac{\varepsilon}{2}$  we finally obtain
\begin{equation}\label{s(1)}
(1)\geq
(\beta-\mu_i)\frac{\varepsilon^2}{4}-C\beta^{\frac{1}{2}+\gamma_1}\qquad\text{for
$\beta$ large enough }.
\end{equation}
Analogously exploiting \eqref{epsilon}, \eqref{51}, \eqref{grad},
 and \eqref{sagmon},
 we discover
that
\begin{equation}\label{s(2)}
(2)\leq C_1,
\end{equation}
where $C_1$ depends on $R$, $S$, and $\|g\|_{W^{1,\infty}}$;
 combining \eqref{s(1)}, \eqref{s(2)}, and recalling \eqref{mu}, we finally
obtain that there exists
 $b_0>\beta_1$ depending only on
 $R$,  $S$,  and $\|g\|_{W^{1,\infty}}$ such that
 \eqref{equivalent} holds true for $\beta\geq b_0$.

Before proceeding let us observe that arguing as for estimate
\eqref{s(2)}, we easily obtain
\begin{equation}\label{sphi}
|\phi_i^x(x,z)|\leq C_2
(\|\nabla
u_{i,\beta}\|_{\infty}+\|\nabla v_i\|_{\infty})\leq C_3 \qquad\text{for
every $(x,z)\in A_i$},
 \end{equation}
where $C_3$ depends only on $R$, $S$,
 and $\|g\|_{W^{1,\infty}}$.
 For $(x,z)\in (\Om_i\times\erre)\cap A_j$ ($i\neq j$),
by the definition of $\phi^x$ and, by \eqref{51},  we have
\begin{equation}\label{preschi}
|\dive_x\phi^x| \leq C_4\left(\|\nabla^2 u_{j,\beta}\|_{\infty}+
 \|\nabla^2
v_j\|_{\infty}+\|\nabla u_{j,\beta}\|^2_{\infty} +
 \|\nabla v_j\|^2_{\infty} +
\|\nabla u_{j,\beta}\|_{\infty}\|\nabla h_{\beta}\|_{\infty}\right) ,
\end{equation}
 where $C_4$ depend only on $R$,
and $S$;
 by using  \eqref{utilde}, \eqref{grad},
\eqref{sagmon}, \eqref{hessian},   \eqref{gradh}, and recalling
that $\gamma_1>\gamma$, we  deduce, from \eqref{preschi}, that
\begin{equation}\label{schi}
|\chi_j|\leq S
C_4\left(C_5\|g\|_{W^{1,\infty}}\beta^{\frac{1}{2}+\gamma}+C_5+K_1\beta_1
+(K\|g\|_{W^{1,\infty}})^2+K_1^2\beta_1 + K\|g\|_{W^{1,\infty}}
K_3\beta^{\frac{1}{2}+\gamma_1}\right)\leq C_6
\beta^{\frac{1}{2}+\gamma_1} ,
\end{equation}
where $C_6$ depends only on $R$, $S$, $\Lambda^{\alpha}(\Gamma)$,
and $\|g\|_{W^{1,\infty}}$ .

Using the definition \eqref{vertical1} for $(x,z)\in
(\Om_1\times\erre)\cap A_2$, we have
\begin{equation}\label{omegai}
\phi^z(x,z)\geq \chi_1-2\|\phi^x\|(\|\nabla u_{2,\beta}\|_{\infty}+
\|\nabla h_{\beta}\|_{\infty})+\phi^z(x,u_{\beta}-h_{\beta}),
\end{equation}
where, by \eqref{verticalparz},
\begin{equation}\label{omegai'}
\phi^z(x,u_{\beta}-h_{\beta})\geq -\beta(u_{\beta}-h_{\beta}-g)^2+
\Psi_1(x,u_{\beta}-h_{\beta}).
\end{equation}

Therefore, for $(x,z)\in
(\Om_1\times\erre)\cap A_2$,  combining
\eqref{omegai} and \eqref{omegai'}, and using \eqref{grad},
\eqref{gradh}, \eqref{spsi},  \eqref{sphi},  and \eqref{schi},
 we obtain
\begin{eqnarray*}
\phi^z(x,z) + \beta(z-g)^2 -\frac{|\phi^x|^2}{4}
&\geq& \beta\left[(z-g)^2-(u_{\beta}-h_{\beta}-g)^2\right]-|\chi_i|-
2\|\phi^x\|_{\infty}(\|\nabla
h_{\beta}\|_{\infty}\\
&& +\|\nabla
u_{2,\beta}\|_{\infty})+
\Psi_1(x, u_{\beta}-h_{\beta})-\frac{\|\phi^x\|_{\infty}^2}{4}\\
&\geq&
\beta[(7/16)^2S^2-(3/16)^2S^2]-C_5\beta^{\frac{1}{2}+\gamma_1}
-C_3(K_3\beta^{\frac{1}{2}+\gamma_1}+
K\|g\|_{W^{1,\infty}})\\
&&-C\beta^{\frac{1}{2}+\gamma_1}-\frac{(C_3)^2}{4},
\end{eqnarray*}
where we used also the fact that
 that $|z-g|\geq
|z-u_{\beta}|-|u_{\beta}-g|\geq S/2-S/16=(7/16)S$ and,
analogously, that
$|u_{\beta}-h_{\beta}-g|\leq
S/16+S/8=(3/16)S$ (see \eqref{vicine}); as
$\frac{1}{2}+\gamma_1<1$, there exists $b_1>0$
depending only on   $R$,  $S$, $\Lambda^{\alpha}(\Gamma)$, and
$\|g\|_{W^{1,\infty}}$ such that
$$
\phi^z(x,z) + \beta(z-g)^2 -\frac{|\phi^x|^2}{4}> 0,
$$
for $\beta\geq b_1$ and for $(x,z)\in
(\Om_1\times\erre)\cap A_2$. Analogously we can prove the existence
of a constant $b_2>0$ depending on the same quantities such that
$$
\phi^z(x,z) + \beta(z-g)^2 -\frac{|\phi^x|^2}{4}> 0,
$$
for $\beta\geq b_2$ and for $(x,z)\in
(\Om_2\times\erre)\cap A_1$.
Arguing exactly in the same way (in fact exploiting the same
estimates), one can finally check  that  there exists $b_3>0$ depending on
$R$,  $S$, $\Lambda^{\alpha}(\Gamma)$, and $\|g\|_{W^{1,\infty}}$,
such that
$$
\phi^z + \beta(z-g)^2 -\frac{|\phi^x|^2}{4}>0
$$
holds true for $(x,z)\in(\Om\times\erre)\setminus
(A_1\cup A_2)$ and $\beta\geq b_3$. If we call $ \beta_2 := \max\{b_0, b_1,b_2, b_3\} $ we
have that for $\beta\geq\beta_2$ condition (c) of Section 1 is satisfied for
almost every $(x,z)$ in
 $\Om\times\erre$ with strict inequality if $z\neq
u(x)$.
 \begin{itemize}
\item $\phi(x,u_{\beta}) = (2\nabla u_{\beta}, |\nabla
u_{\beta}|^2-\beta(u_{\beta}-g)^2)$ {\it everywhere in}
$\Om\setminus\Gamma$.
\end{itemize}
 Condition (d) of Section 1 is trivially satisfied, as one
 can see directly from the definition of $\phi$.
 \begin{itemize}
\item{$\displaystyle\int_{u_{2,\beta}(x)}^{u_{1,\beta}(x)}\phi^x(x,t)\, dt =
\nu_{u_{\beta}}=-\nu$ ${\cal H}^{n-1}$-{\it a.e. on} $\Gamma$.}
\end{itemize}
By direct
computation, for $x\in\Gamma$, we have
\begin{equation}\label{star}
\int_{u_{2,\beta}}^{u_{1,\beta}}\phi^x(x,z)\, dz =(h_{\beta})^2
\frac{\nabla v_2}{v_2} -  (h_{\beta})^2\frac{\nabla v_1}{v_1}.
\end{equation}
Using \eqref{normale}, \eqref{tildeh}, and the fact that
$v_i\equiv 1$ on $\Gamma$, we obtain
$$
\int_{u_{2,\beta}}^{u_{1,\beta}}\phi^x(x,z)\, dz=
+\frac{1}{2}\frac{\nabla v_2}{|\nabla v_2|}-
\frac{1}{2}\frac{\nabla v_1}{|\nabla v_1|} = -\nu,$$ so that
condition (e) of Section 1 is satisfied.
\begin{itemize}
\item{$\left|\displaystyle\int_{t_1}^{t^2}\phi^x(x,z)\, dz\right|\leq1$ {\it
for every}
 $t_1$, $t_2\in\erre$ {\it and for every} $x\in\Om$. }
\end{itemize}
It is convenient to introduce the following notation: for every $x\in\Om$ and
for every $s$, $t\in\erre$, we set
$$
I(x,[s,t]):= \int_s^t\phi^x(x,z)\, dz,
$$
where, with a slight abuse of notation, $[s,t]$ stands
for the interval $[s\land t, s\lor t]$ positively oriented
if $s\leq t$, negatively oriented otherwise.
We define
$$
d_{\beta}(\pi(x)):=\frac{\tilde
h(\pi(x))-\varepsilon}{\beta^{\frac{1}{2}+\gamma_1}}.
$$

 If $|d(x)|> d_{\beta}(\pi(x))$, recalling that, by
 definition, $h_{\beta}(x)=\varepsilon$
we have
 \begin{eqnarray}
|I(x,t_1,t_2)|&\leq & \int_{ u_{1,\beta}-\varepsilon}^{
u_{1,\beta}+\varepsilon}\left(\hbox{$2\|\nabla
u_{1,\beta}\|_{\infty}+\frac{16}{\varepsilon} \|\nabla
u_{1,\beta}\|_{\infty}\left(
u_{1,\beta}-\frac{\varepsilon}{2}-z\right)^+ +4{|
u_{1,\beta}-z|}\|\nabla
v_1\|_{\infty} $} \right )\, dz\nonumber\\
&&+\int_{ u_{2,\beta}-\varepsilon}^{
u_{2,\beta}+\varepsilon}\left(\hbox{$2\|\nabla
u_{2,\beta}\|_{\infty}+\frac{16}{\varepsilon} \|\nabla
u_{2,\beta}\|_{\infty}\left(
-u_{2,\beta}-\frac{\varepsilon}{2}+z\right)^+ +4{|
u_{2,\beta}-z|}\|\nabla
v_2\|_{\infty} $} \right )\, dz\nonumber\\
 &\leq& 6\varepsilon\|\nabla  u_{1,\beta}\|_{\infty} +4
\varepsilon^2\|\nabla v_1\|_{\infty}+6\varepsilon\|\nabla
u_{2,\beta}\|_{\infty} +4 \varepsilon^2\|\nabla v_2\|_{\infty}
\leq \frac{1}{2},\label{sepsilon}
\end{eqnarray}
where the last inequality is due to \eqref{epsilon}, therefore
condition (f) is satisfied.

Let us consider now the case of a point $x$ where $|d(x)|\leq
d_{\beta}(\pi(x))$.
For $x\in \Om_i\cup\Gamma$ we set
$$n(x):= - \frac{\nabla
v_1}{|\nabla v_1|}=\frac{\nabla
v_2}{|\nabla v_2|};$$
 note that $n(x) = \nu_{u_{\beta}}(x)$ for every
$x\in\Gamma$. Given any vector valued function
$\xi:\Om\to\erre^n$, we call $\xi^{\perp}$ and $\xi^{\parallel}$ the vector
valued functions such that $\xi^{\perp}(x)$ and $\xi^{\parallel}(x)$ are equal
to the projections of
$\xi(x)$ on the orthogonal space and
on the space generated by $n(x)$, respectively. We denote the open unit
sphere of $\erre^n$  centred at the origin by $B$ and the open ball of
$\erre^n$
 centred at the point $-rn(x)$ with radius $r$, by $A(x,r)$.
Finally, for $x\in\Om$ and $t\in\erre$ we introduce the following
 vector
$$ b_i(x,t):=(-1)^{i}(2(t- u_{i,\beta})-
j_i(x,t))\left(\nabla u_{i,\beta}\right)^{\parallel},
$$
where $j_i$ is defined by
$$
j_i(x,t):=\frac{16}{h_{\beta}}\int_{ u_{i,\beta}}^t
\hbox{$\left( (-1)^i(u_{i,\beta}-z)-\frac{h_{\beta}}{2}\right)^+
$}\, dz.
$$
CLAIM 1. There exists a positive constant $c_0>0$, depending on
$R$, $S$, $\Lambda^{\alpha}(\Gamma)$, and $\|g\|_{W^{1,\infty}}$,
with the property that for  every $x\in\Om$ such that
$|d(x)|\leq d_{\beta}(\pi(x)) $, for every $t\in\erre$ such that $|t-
u_{i,\beta}(x)|\leq h_{\beta}(x)$, and for $\beta\geq c_0$, we
have
\begin{equation}\label{171}
(-1)^{i+1}I(x, [u_{i,\beta},t])+ b_i(x,t) \in A(x,1/3).
\end{equation}

\noindent
A straightforward computation gives
\begin{multline}
(-1)^{i+1}I(x, [u_{i,\beta}(x),t])+ b_i(x,t)=
2(-1)^{i+1}\nabla  u_{i,\beta} (t- u_{i,\beta})+
(-1)^ij_i(x,t)\nabla  u_{i,\beta}+
 (-1)^{i+1}\frac{\nabla v_i}{v_i}(t-u_{i,\beta})^2\\
+2(-1)^{i}(t-
u_{i,\beta})\left(\nabla  u_{i,\beta}\right)^{\parallel}+
(-1)^{i+1}j_i(x,t)\left(\nabla  u_{i,\beta}\right)^{\parallel}=\\
 = (-1)^{i+1}(2(t-u_{i,\beta})-j_i(x,t))(\nabla  u_{i,\beta})^{\perp}
 -\frac{|\nabla v_i|}{v_i}(t-u_{i,\beta})^2n(x) \nonumber
\end{multline}
and so the claim is equivalent to prove that
$$
(2-j_i(x,t)(t- u_{i,\beta})^{-1})^2[(\nabla
u_{i,\beta})^{\perp}]^2(t- u_{i,\beta})^2+\frac{|\nabla
v_i|^2}{v_i^2} (t-u_{i,\beta})^4 - \frac{2}{3}\frac{|\nabla
v_i|}{v_i}(t-u_{i,\beta})^2<0;
$$
as $0\leq 2-j_i(x,t)(t-u_{i,\beta})^{-1}\leq 2$
everywhere, it is sufficient to  prove that
\begin{equation}\label{tesieq}
(*):=4\left|(\nabla u_{i,\beta})^{\perp}\right|^2+ h_{\beta}^2
\frac{|\nabla v_i|^2}{ v_i^2} - \frac{2}{3}\frac{|\nabla
v_i|}{v_i}<0.
\end{equation}
Since, by \eqref{tildeh}, $h^2_{\beta}\frac{|\nabla
v_i|}{v_i}=\frac{1}{2}$ for $x\in\Gamma$, we can estimate
\begin{equation}\label{181}
(*)= 4[(\nabla u_{i,\beta})^{\perp}]^2-\frac{1}{6}\frac{|\nabla
v_i|}{v_i}<0 \qquad\text{on $\Gamma$},
\end{equation}
where the last inequality follows from \eqref{v12}. In the
following we denote by $\partial_{|d|}$ the differential operator
$$
\partial_{|d|}f(x)=\nabla f(x)\cdot\nabla |d(x)|,
$$
defined for $x\in(\Gamma)_{\frac{R}{2}}\setminus\Gamma$; noting
that, by the estimates \eqref{sagmon}, we have
$$
\left|\partial_{|d|}\frac{|\nabla v_i |^2}{v_i^2}\right|\leq
C\qquad \left|\partial_{|d|}\frac{|\nabla v_i |}{v_i}\right|\leq
C,
$$
where $C$ depends only on $R$, $S$, and $\|g\|_{W^{1,\infty}}$,
and using \eqref{Dtilde}, \eqref{epsilon}, \eqref{hessian},
\eqref{utilde}, and \eqref{gradh}, one sees that
\begin{eqnarray*}
\partial_{|d|}((*)) & = & 8(\nabla u_{i,\beta})^{\perp}\cdot
\partial_{|d|}(\nabla  u_{i,\beta})^{\perp}+2\frac{|\nabla v_i|^2}{ v_i^2}
h_{\beta}\partial_{|d|}h_{\beta}+h^2_{\beta}\partial_{|d|}\frac{|\nabla
v_i |^2}{v_i^2}
-\frac{2}{3}\partial_{|d|}\frac{|\nabla v_i |}{v_i}\\
&\leq& 8cKK_2\beta^{\frac{1}{2}+\gamma}\|g\|_{W^{1,\infty}}^2+C_1 -
\frac{\varepsilon}{2} K_3\beta^{\frac{1}{2}+\gamma_1}  +S^2C+C;
\end{eqnarray*}
as $\gamma_1>\gamma$ and since all the constants appearing in the
last inequality depend only on $R$, $S$,
$\Lambda^{\alpha}(\Gamma)$, and $\|g\|_{W^{1,\infty}}$,   it is
clear that there exists $c_0>0$ depending on the same quantities
such that $
\partial_{|d|}((*))<0,
$ for $x\not\in\Gamma$ such that $|d(x)|\leq d_{\beta}(\pi(x))$ and for
$\beta\geq c_0$. Therefore, taking into account  \eqref{181},
\eqref{tesieq} follows immediately: Claim 1  is proved.

CLAIM 2. There exists a positive constant $c_1$, depending only on
 $R $, $S$, $\Lambda^{\alpha}(\Gamma)$, and
$\|g\|_{W^{1,\infty}}$, such that
for every $x\in\Om$, $t_1$, $t_2\in\erre$, with
$|d(x)|\leq  d_{\beta}(\pi(x))$, $|t_1-u_{1,\beta}|\leq h_{\beta}$,
$|t_2-u_{2,\beta}|\leq h_{\beta}$, and for every $\beta\geq c_1$, we have
\begin{equation}\label{claim2}
I(x,[u_{2,\beta},
u_{1,\beta}])-b_1(x,t_1)-b_2(x,t_2)=b(x,t_1,t_2)n(x),
\end{equation}
with $b(x,t_1,t_2)<1$.

First of all observe that for every $x\in\Om$ $I(x,
u_{2,\beta}, u_{1,\beta})$ is a vector parallel to $n(x)$,
by \eqref{star};
it is also clear that
\begin{eqnarray*}
|I(x,[u_{2,\beta},
u_{1,\beta}])-b_1(x,t)-b_2(x,t)|&\leq& |I(x,[
u_{2,\beta}, u_{1,\beta}])| + \left|\left(\nabla
u_{1,\beta}\right)^{\parallel}\right| |(2h_{\beta} + j_1(x,h_{\beta}))\\
&&+\left|\left(\nabla u_{2,\beta}\right)^{\parallel}\right|(2h_{\beta} +
j_2(x,h_{\beta}))\\
&\leq&|I(x,[u_{2,\beta},u_{1,\beta}])|+ 4h_{\beta}
\left(\left|\left(\nabla
u_{1,\beta}\right)^{\parallel}\right|+\left|\left(\nabla
u_{2,\beta}\right)^{\parallel}\right|\right)\\
&&=:m_{\beta}(x);
\end{eqnarray*}
therefore it is sufficient to prove that $m_{\beta}(x)<1$
for $|d(x)|\leq d_{\beta}(\pi(x))$, if $\beta$ is large enough. Since
$m_{\beta}(x)=|I(x, u_{2,\beta}, u_{1,\beta})|=1$ for
every $x\in\Gamma$, it will be enough to show that
$\partial_{|d|}m_{\beta}(x)<0$ for $x$ such
that $|d(x)|\leq d_{\beta}(\pi(x))$. We don't enter all the details, indeed
arguing as above, that is using \eqref{grad}, \eqref{hessian},
\eqref{sagmon}, and \eqref{gradh}, one easily sees that the
derivative of $h_{\beta}$ which is negative and of the same order
as $\beta^{\frac{1}{2}+\gamma_1}$, dominates the other terms and
so there exists a positive constant $c_1>0$ depending on $R$,
$S$, $\Lambda^{\alpha}(\Gamma)$, and $\|g\|_{W^{1,\infty}}$, such
that $\partial_{|d|}m_{\beta}(x)<0$ for $\beta\geq c_1$: Claim 2
is proved.

We set $\beta_3= \max\{c_0,c_1\}$ and we are going to prove that
condition (f) of Section 1 is satisfied for
$\beta\geq\beta_3$. We will check the condition only in
$\Om_1\times\erre$: for $\Om_2\times\erre$ the argument would be
analogous. Let $x\in\Om_1$ and $t_2<t_1$ two real numbers such
that $|t_2- u_{2,\beta}(x)|\leq h_{\beta}(x)$ and
$|t_1-u_{\beta}(x)|\leq h_{\beta}(x)$; first of all it is easy
to see, by explicit computation, that
\begin{equation}\label{perfortuna}
I(x,[t_2,t_1])\cdot n(x)\geq0;
\end{equation}
recalling that, by Claim 1,
$$
I(x,[u_{\beta},t_1])+b_1(x,t_1)\in
A(x,1/3)\qquad\text{and}\qquad I(x,[t_2,
u_{2,\beta}])+b_2(x,t)  \in A(x,1/3),
$$
we have
\begin{eqnarray*}
I(x,[t_2,t_1])  &=& I(x,[
u_{2,\beta},u_{\beta}])-b_1(x,t_1)-b_2(x,t_2) + I(x,[t_2,
u_{2,\beta}])+b_2(x,t_2)\\
&& + I(x,[u_{\beta},t_1])+b_1(x,t_1)\in I(x,[
u_{2,\beta},u_{\beta}])-b_1(x,t_1)-b_2(x,t_2)+ 2A(x,1/3),
\end{eqnarray*}
therefore,  taking
into account \eqref{perfortuna},
\begin{equation}\label{cisiamo}
I(x,[t_2,t_1])\in (I(x,
[u_{2,\beta},u_{\beta}])-b_1(x,t_1)-b_2(x,t_2) + A(x,2/3))\cap
H^+,
\end{equation}
where $H^+$ is the half-space $\{\xi\in\erre^n: \xi\cdot n(x)\geq
0\}$. By elementary geometry it is easy to see that
$(bn(x)+A(x,r))\cap H^+\subset B$ for $b<1$ and for
$r\in(0,1)$, and hence, invoking  Claim 2, we get
\begin{eqnarray}
I(x,[t_1,t_2])&\in&\left(I(x,[
u_{2,\beta},u_{\beta}])-b_1(x,t_1)-b_2(x,t_2) +
A(x,2/3)\right)\cap H^+\nonumber\\
&&= (b(x,t_1,t_2)n(x)+
A(x,2/3))\cap H^+\subset B.\label{perora}
\end{eqnarray}
If $(x,t_1)$ and $(x,t_2)$ belong to $A_i$ it is easy to see, by
explicitly computing the integral, that
\begin{equation}\label{lili}
|I(x,[t_1,t_2])| \leq h_{\beta}^2(x)\frac{|\nabla
v_i|}{v_i}+\frac{25}{8}h_{\beta}|\nabla u_{i,\beta}|<
1-\frac{25}{32}\frac{1}{\sqrt3}+\frac{25}{32}\frac{1}{\sqrt3}=1,
\end{equation}
where the last inequality follows from \eqref{Dtilde}, \eqref{v12bis}, and
\eqref{tildeh}
(we recall that for $\beta$ large
enough $d_{\beta}(\pi(x))\leq D$, for every $x$, being $D$ the constant
introduced in \eqref{Dtilde}).

We now consider the general case. Let $x\in\Om_1$, $t_1$, $t_2\in\erre$ with
$t_1<t_2$; since $\phi^x$ vanishes out of the regions $A_1$ and $A_2$, we have
$$
I(x,[t_1,t_2])=I(x,[t_1,t_2]\cap[u_{2,\beta}-h_{\beta},
u_{2,\beta}+h_{\beta}])+I(x,[t_1,t_2]\cap[u_{\beta}-h_{\beta},
u_{\beta}+h_{\beta}]);
$$
by \eqref{lili}, each integral in the expression above has modulus
less than 1, so that if one of the two is vanishing condition (f) is verified.
If both are non-vanishing, then
$$
[t_1,t_2]\cap[u_{2,\beta}-h_{\beta},u_{\beta}+h_{\beta}]=[s_1,s_2],
$$
with $|s_1-u_{2,\beta}|\leq h_{\beta}$ and $|s_2-u_{\beta}|\leq
h_{\beta}$, so that, again taking into account the fact that  $\phi^x$
vanishes out of the regions $A_1$ and $A_2$,
 $$
|I(x,[t_1,t_2])|=
|I(x,[t_1,t_2]\cap[u_{2,\beta}-h_{\beta},u_{\beta}+h_{\beta}])|=
|I(x,[s_1,s_2])|<1,
$$
where the last inequality follows from \eqref{perora}: condition (f)
of Section 1 is
proved.

Since, by construction, $\phi$ has vanishing normal component on
$\partial\Om\times\erre$, if we set $\overline
\beta:=\max\{\beta_1,\beta_2,\beta_3\}$ we have that conditions of
Section 1 are all satisfied for
$\beta\geq\overline\beta$: the theorem is proved. \qed

A similar result holds true also if $\Gamma$ is made up of several
connected components, as the following theorem states: we omit the
proof, since it is essentially the same as the previous one.

\begin{theorem}\label{cor1}
Let $\Om$ as above and let $\Om_1,\dots,\Om_k$ a family of  open
disjoint subsets belonging of class
$C^{2,\alpha}$ and let $R>0$ such that
$\Om_i\in {\cal U}_R(\Om)$ for $i=1,\dots,k$ and
$\dist(\Om_i,\Om_j)\geq R$ for every $i\neq j$.
 Set
$\Gamma:=\partial\Om_1\cup\dots\cup\partial\Om_k$.  Then
for every function $g$ belonging
$W^{1,\infty}(\Om\setminus\Gamma)$, discontinuous along $\Gamma$
(i.e. $S_g=\Gamma$) and such that $g^+(x)-g^-(x)>S>0$ for every
$x\in\Gamma$, there exists $\beta_0>0$ depending on $R$, $S$,
 $\Lambda^{\alpha}(\Gamma)$ (see \eqref{normaholder}), and
$\|g\|_{W^{1,\infty}}$, such that for $\beta\geq\beta_0$ the
solution $u_{\beta}$ of \eqref{ubeta'} is discontinuous along
$\Gamma$ ($S_{u_{\beta}}=\Gamma$) and it is the unique absolute
minimizer of $F_{\beta,g}$ over $SBV(\Om)$.
\end{theorem}

\begin{remark}
We remark that refining a little the construction, it is possible
to improve the result of Theorem \ref{calibration1} as follows:

\noindent {\it there exist $\delta^*>0$ and $\beta_0>0$ such that,
for every $\beta\geq\beta_0$ and  for every $g\in
W^{1,\infty}(\Om\setminus\Gamma)$, with $\|g\|_{W^{1,\infty}}\leq
\beta^{\delta^*}$ and such that $\inf_{\Gamma}(g^+-g^-)>S$, the
solution $u_{\beta}$ of \eqref{ubeta'} is the unique absolute
minimizer of $F_{\beta,g}$ over $SBV(\Om)$. }

\noindent The main difficulty comes from the fact that instead of
\eqref{grad} we have the weaker estimate
$$
\|\nabla u_{\beta}\|_{\infty}\leq K \beta^{\delta^*}.
$$
Such a difficulty can be overcome replacing, in the construction
above, $v_1$ and $v_2$ by  $v_{1,\beta}$ and
$v_{2,\beta}$ defined as
$$
v_{1,\beta}(x)=
\begin{cases}
z_{1,c\beta^{4\delta^*}}(x) & \text{if $x\in\overline{\Om_1}$}\\
2-z_{2,c\beta^{4\delta^*}}(x) &\text{if $x\in\Om_2$}
\end{cases}
$$
and
$$
v_{2,\beta}(x)=
\begin{cases}
z_{2,c\beta^{4\delta^*}}(x) & \text{if $x\in\overline{\Om_2}$}\\
2-z_{1,c\beta^{4\delta^*}}(x) &\text{if $x\in\Om_1$},
\end{cases}
$$
where $z_{1,c\beta^{4\delta^*}}$ and $z_{2,c\beta^{4\delta^*}}$
are the two functions constructed in Lemma \ref{derivnorm} with
$\lambda=c\beta^{4\delta^*}$. One can check that if $\delta^*$ is
sufficiently small and  $c$  sufficiently large, all the
conditions of Section 1 are still satisfied for
$\beta$ large enough.
\end{remark}

\subsection{The two-dimensional case}

As stated in the Introduction, in dimension two  we are able to
treat the case of $\Om$ with
piecewise smooth boundary (curvilinear polygon) and of $\Gamma$
touching (orthogonally) $\partial\Om$.
\begin{lemma}\label{derivnorm2}
Let $\Om$, ${\cal S}$, and $\Gamma$ be as in Proposition
\ref{BOOO} and denote by  $\Om_1$, $\Om_2$ the two connected
components of $\Om\setminus\Gamma$. Then for every $\delta>0$
there exist two positive constants $c$ and $\beta_0$ depending on
$\Gamma$ and $\delta$ (and $\Om$ of course) such that, for
$\beta\geq \beta_0$, we can find two functions
$z_{1,\beta}:\Om_1\to\erre$ and $z_{2,\beta}:\Om_2\to\erre$ of
class $W^{2,\infty}$ with the following properties:
\begin{description}
\item{i)}$\frac{1}{2}\leq z_{i,\beta}\leq 1$ in $\Om_i$, for $i=1,2$ and
$z_{i,\beta}\equiv \frac{1}{2}$ in
$\Om\setminus(\Gamma)_{\delta}$;
\item{ii)}$\Delta z_{i,\beta}\leq c\beta z_{i,\beta}$ in $\Om_i$, for
$=1,2$;
\item{iii)} $z_{1,\beta}(x)=z_{2,\beta}(x)=1$ and
$\partial_{\nu} z_{1,\beta}(x)=-\partial_{\nu}
z_{2,\beta}(x)\geq\sqrt\beta$ for every $x\in\Gamma$;
\item{iv)} $\|\nabla z_{i,\beta}\|_{\infty}\leq c\sqrt{\beta}$ and
$\|\nabla^2 z_{i,\beta}\|_{\infty}\leq c{\beta}$.
\end{description}
\end{lemma}
{\sc Proof.} Let us denote by $x_1$ and $x_2$ the two intersection
points of $\Gamma$ with $\partial\Om$. If we are able to find a
function $\tilde{d}$ belonging to
$W^{2,\infty}((\Gamma)_{\delta'}\cap\Om)$ (for a suitable
$\delta'<\dist({\cal S},\Gamma)$) such that $\tilde d$ is
vanishing on $\Gamma$, positive in $\Om_2\cap (\Gamma)_{\delta'}$,
negative in  $\Om_1\cap (\Gamma)_{\delta'}$, satisfying
$\partial_{\nu}\tilde d =0$ on
$\partial\Om\cap\overline{\Gamma^{\delta'}}$ and
$\partial_{\nu}\tilde d \neq0$ on $\Gamma$, we are done: indeed we
can proceed exactly as in Lemma \ref{derivnorm} using $\tilde d$
in place of $d$.  We briefly describe a possible construction: as
in Proposition \ref{BOOO} we can find a neighbourhood $U_i$ of
$x_i$ ($i=1,2$) and a $C^{1,1}$ function $\psi_i$ vanishing on
$\Gamma\cap U_i$, positive in $\Om_2\cap U_i$, negative in
$\Om_1\cap U_i$ and such  that $\partial_{\nu}\psi_i =0$ on
$\partial\Om\cap U_i$ and $\partial_{\nu}\psi_i \neq 0$ in
$\Gamma\cap U_i$. Now we can define $\tilde d:= \theta_1
\psi_1+\theta d+\theta_2 \psi_2$, where $\theta_1$, $\theta_2$,
and $\theta_3$ are suitable positive cut-off functions such that
$\theta_1+\theta_2+\theta_3\equiv 1$, while $d$ is the usual
signed distance function from $\Gamma$, positive in $\Om_2$ and
negative in $\Om_1$ (it is well defined in $\Gamma^{\delta'}$ if
$\delta'$ is small enough). \qed

\begin{theorem}\label{calibration2}
Let $\Om$, $\Om_1$, $\Om_2$, and $\Gamma$ as in the previous Lemma
and let  $g$ be a function in  $W^{1,\infty}(\Om\setminus\Gamma)$,
discontinuous along $\Gamma$ (i.e. $S_g=\Gamma$) and such that
$g^+(x)-g^-(x)>S>0$ for every $x\in\Gamma$. Then there exists
$\beta_0>0$ depending on $\Gamma$, $S$, and
$\|g\|_{W^{1,\infty}}$, such that for $\beta\geq\beta_0$ the
solution $u_{\beta}$ of \eqref{ubeta'} is discontinuous along
$\Gamma$ ($S_{u_{\beta}}=\Gamma$) and it is the unique absolute
minimizer of $F_{\beta,g}$ over $SBV(\Om)$.
\end{theorem}
{\sc Proof.} As above, let us denote by ${\cal S}$ the set of the
singular points of $\partial\Om$. If $\Om$ is regular (i.e. ${\cal
S}=\emptyset$) we can recycle exactly the same construction of
Theorem \ref{calibration1}. If ${\cal S}\neq\emptyset$, an
additional difficulty is due to the fact that we are not able to
prove that $\|\nabla u_{\beta}\|_{L^{\infty}(\Om)}\leq C$ with $C$
independent of $\beta$. Since we can perform such an estimate only
in a neighbourhood of $\Gamma$ which does not intersect ${\cal
S}$, the idea will be to keep the construction of Theorem
\ref{calibration1} in that neighbourhood  and to suitably modify
it near the singular points in order to exploit estimate
\eqref{gradirreg}.

Denote by $\gamma_1$ and $\gamma_2$ the two curvilinear edges of
$\Om$ containing the intersection points of $\Gamma$ with
$\partial\Om$ and choose $\delta>0$ so small that
$(\Gamma)_{\delta}\cap {\cal S}=\emptyset$,
$(\Gamma)_{\delta}\cap\partial\Om=(\Gamma)_{\delta}\cap(\gamma_1\cup\gamma_2)$,
and  $d$ and $\pi$ are well defined and smooth (according to Lemma
\ref{-2}) in that neighbourhood.

Let us choose $\beta'>0$ and $G>0$
such that, for $\beta\geq\beta'$,
\begin{equation}\label{vicine2}
\|u_{\beta}-g\|_{L^{\infty}(\Om)}\leq \frac{S}{16}
\quad\text{and}\quad\sqrt\beta\|u_{\beta}-g\|_{L^{\infty}(\Om)}\leq
G\|g\|_{W^{1,\infty}(\Om)}\quad i = 1,2,:
\end{equation}
this is possible by virtue of Proposition \ref{prop2}.

Again it is convenient to extend the restriction of $u_{\beta}$
to $\Om_i$
($i=1,2$) to a $C^{1,1}$ function $ u_{i,\beta}$ defined in
the whole $\Om$, in such a way that
\begin{equation}\label{utilde2}
u_{i,\beta}(x) = u_{\beta}(x) \ \text{in $
{\Om}_i$,}\quad \|
u_{i,\beta}\|_{W^{2,\infty}(\Om)} \leq c
\|u_{\beta}\|_{W^{2,\infty}((\Gamma)_{\delta}\cap\Om)},
\quad\text{and}\quad u_{1,\beta}-u_{2,\beta} \geq
\frac{3}{4}S\quad \text{everywhere},
\end{equation}
where $c$  is a positive constant independent of $\beta$. We
require also that
$$
\partial_{\nu} u_{i,\beta} = 0\qquad\text{on
$\partial\Om$}.
$$
 By \eqref{scasoirreg} and \eqref{utilde2}, and by \eqref{gradirreg}, we can
state the existence of two positive constants $K$ and $\beta''$
depending only on $\Gamma$,
 such that
\begin{equation}\label{gradu2}
\|\nabla u_{i,\beta}\|_{L^{\infty}(\Om_{\io}\cup
(\Gamma)_{\delta}\cap\Om)}\leq
K\|g\|_{W^{1,\infty}(\Om)}\quad\text{for}\quad
i=1,2,\qquad\text{and}\qquad \|\nabla
u_{\beta}\|_{L^{\infty}(\Om)}\leq\beta^{\frac{1}{4}}
K\|g\|_{W^{1,\infty}(\Om)}
\end{equation}
for every $\beta\geq \beta''$ (above and in the sequel,
$\io$ denotes the comlement of $i$, i.e., $\io$ is such that
${i,\io}=\{1,2\}$).

Let $\beta'''>0$ satisfying
\begin{equation}\label{beta'''2}
\frac{1}{6}\sqrt{\beta'''}=\max\left\{4(K\|g\|_{W^{1,\infty}})^2,
{64}/{S^2},\beta',\beta'',\beta_0\right\}+1,
\end{equation}
where $\beta_0$ is the constant appearing in  Lemma
\ref{derivnorm2} and  $z_{1,\beta'''}$, and let $z_{2,\beta'''}$
be the two functions constructed in Lemma \ref{derivnorm2} with
$\lambda=\beta'''$. We denote by $v_{1}$, $v_{2}$ the functions
defined as follows:
$$
v_{1}(x)=
\begin{cases}
z_{1,\beta'''}(x) & \text{if $x\in\overline{\Om_1}$}\\
2-z_{2,\beta'''}(x) &\text{if $x\in\Om_2$}
\end{cases}
$$
$$
v_{2}(x)=
\begin{cases}
z_{2,\beta'''}(x) & \text{if $x\in\overline{\Om_2}$}\\
2-z_{1,\beta'''}(x) &\text{if $x\in\Om_1$},
\end{cases}
$$
and we choose $0<D<\delta$ in such a way that
$$
|\nabla v_i(x)|\geq\frac{1}{2},
\quad \tilde h^2(\pi(x))\frac{|\nabla v_i|}{v_i}\leq 1-\frac{25}{32}
\frac{1}{\sqrt3},
\qquad\text{if $|d(x)|\leq D$,
$i=1,2$,}
$$
where
\begin{equation}\label{tildeh2}
\tilde h(x)=\frac{1}{\sqrt{2}}|\nabla v_1(x)|
^{-\frac{1}{2}}=\frac{1}{\sqrt{2}}|\nabla v_2(x)|
^{-\frac{1}{2}}\qquad\forall x\in\Gamma.
\end{equation}
 Then we
choose $\varepsilon\in (0,1)$ in such a way that
\begin{equation}\label{epsilon2}
 12\varepsilon\|\nabla\tilde
u_{i,\beta}\|_{L^{\infty}(\Om_{\io}\cup(\Gamma)_{\delta})\cap\Om}+4\varepsilon^2
\|\nabla v_i\|_{L^{\infty}(\Om)}<\frac{1}{4}
\quad \text{for $i=1,2$
and $\beta\geq
 \beta'''$}.
\end{equation}

 Let $\gamma$ be a fixed constant belonging to $(0,\frac{1}{2}\land\alpha)$:
by Proposition \ref{BOOO}, we can find two positive constants
$\beta^{\text{\i v}}$ and $K_2$  such that
\begin{equation}\label{hessian2}
\|\nabla^2
u_{\beta}\|_{L^{\infty}\left((\Gamma)_{\delta}\cap\Om\right)}\leq
K_2 \beta^{\frac{1}{2}+\gamma}\|g\|_{W^{1,\infty}(\Om)},
 \end{equation}
for every $\beta\geq\beta^{\text{\i v}}$.

Now we can
define, for $\beta>0$,
$$
 h_{\beta}(x) =
\begin{cases}
\left(\tilde
h(\pi(x))-\beta^{\frac{1}{2}+\gamma_1}|d(x)|\right)\lor
\varepsilon & \text{if $|d(x)|\leq \frac{D}{2}$}\\
 f_{\beta}(|d(x)|)&\text{if $|d(x)|>\frac{D}{2}$},
\end{cases}
$$
where $\gamma_1$ is a fixed constant belonging to
$(\gamma,\frac{1}{2})$ and $f_{\beta}:[\delta,+\infty)\to \erre$ is the
continuous function satisfying
\begin{equation}\label{fbeta}
f_{\beta}\left(\frac{D}{2}\right)=\varepsilon\qquad
f_{\beta}(t)\equiv s_{\beta}
:=\left({\beta^{\frac{1}{4}}8K\|g\|_{W^{1,\infty}(\Om)}}\right)^{-1}
\quad\text{for
$t\geq D$} \qquad \text{$f_{\beta}$ is affine in
$\left[\frac{D}{2},D\right]$}.
\end{equation}
It is easy to see that there exists
$\beta^{\rm v}>0$  such that   $h_{\beta}$ is continuous
(in fact Lipschitz) for $\beta>\beta^{\rm v}$.

Finally we introduce a new function $\hat u_{i,\beta}$  which is a
modification of  $u_{i,\beta}$ in the region where we cannot perform a
uniform  control of the $L^{\infty}$-norm of its gradient; such a function must
satisfy, for $i=1,2$:
\begin{equation}\label{ubetamod1}
\hat u_{i,\beta}(x)= u_{i,\beta}(x)\quad\text{for
$x\in\Om_{\io}\cup(\Gamma)_{\frac{D}{2}}\cap\Om$}\qquad\text{and}\qquad
 \hat u_{i,\beta}(x)=g(x)\quad\text{for
$x\in\Om_i\setminus(\Gamma)_{{D}}$,}
\end{equation}
\begin{equation}\label{ubetamod2}
\|\nabla \hat u_{i,\beta} \|_{L^{\infty}(\Om)}\leq c (\|\nabla
u_{\beta}\|_{L^{\infty}((\Gamma)_{D})}\lor \|\nabla
g\|_{L^{\infty}(\Om)})\qquad\text{and}\qquad \|\hat
u_{i,\beta}-g\|_{L^{\infty}(\Om)} \leq
\|u_{\beta}-g\|_{L^{\infty}(\Om)},
\end{equation}
where $c>0$ is independent of $\beta$: a possible construction is given by
$$\hat u_{i,\beta}(x)=\theta\left((-1)^id(x)\right)\tilde
u_{i,\beta}+\left[1-\theta\left((-1)^i d(x)\right)\right]g(x),$$
where $\theta$ is a smooth positive function such  that
$\theta(t)= 1$ for $t\leq D/2$ and $\theta(t)=0$ for $t\geq D$.
Now for $\beta\geq \beta_1 :=
\max\{\beta'',\beta''',\beta^{\text{\i v}},\beta^{\rm v  }\} $ we
consider the sets
\begin{equation}\label{Ai2} A_i := \{(x,z)\in\Om\times\erre:\
\hat u_{i,\beta}(x)-h_{\beta}(x)\leq z\leq\hat u_{i,\beta}(x) +
h_{\beta}(x)\}, \quad i = 1,2;
\end{equation}
setting
$$\hat
h_{\beta}(x):=\left[1+(2/D)(|d(x)|-D/2)^+\right]h_{\beta}(x),$$
we can define
 $$
\phi^x(x,z): =
\begin{cases}
2\nabla u_{i,\beta}-2\frac{u_{i,\beta}-z}{v_i}\nabla v_i -\frac{16}{\hat
h_{\beta}}\left((-1)^i(z-u_{i,\beta})-\frac{\hat
h_{\beta}}{2}\right)^+\nabla u_{i,\beta}&\text{if $(x,z)\in A_i$},\\
0 & \text{otherwise,}
\end{cases}
$$
and
$$
\phi^z|_{A_i\cap(\overline \Om_i\times\erre)}:=\left|\nabla
u_{\beta}-\frac{u_{\beta}-z}{v^i}\nabla v^i\right|^2-
 \beta(z-g)^2  +(\beta-\mu_i)(u_{\beta}-z)^2
+ \Psi_{i},
$$
where the functions $\Psi_{i}$ and $\mu_i$ are defined exactly as in
the proof of Theorem \ref{calibration1}. At this point, as in
the proof of Theorem \ref{calibration1}, the vertival component
$\phi^z$ can be extended to the whole $\Om\times\erre$ in order to
satisfy conditions (a) and (b) of Section 1 (we do not rewrite
the explicit expression).
First of all observe that in
$\left((\Gamma)_{\frac{D}{2}}\cap\Omega\right)\times\erre$ the
definition of $\phi$ is the same as in Theorem \ref{calibration1}, then,
we can state the existence of a constant $\beta_0'>0$ depending on
$\Gamma$, $S$, and $\|g\|_{W^{1,\infty}}$ such that $\phi$
satisfies (a), (b), (c), (d),  (e), (f), and (g) of Section 1 in
$\left((\Gamma)_{\frac{D}{2}}\cap\Omega\right)\times\erre$.

From now on we focus our attention on what
happens in
$\left(\Om_i\setminus(\Gamma)_{\frac{D}{2}}\right)\times\erre.$

Concerning (d), we have only to check that for $\beta$ large enough
the graph of $u_{i,\beta}$ belongs to $A_i$, but this follows from the
fact that, by \eqref{fbeta}, $A_i$ contains the $s_{\beta}$-neighbourhood
of the graph of $\hat u_{i,\beta}$, where $s_{\beta}$ is of order
$\beta^{-\frac{1}{4}}$, and from the fact that, by \eqref{vicine2} and
\eqref{ubetamod2}, it holds
$$
\|u_{i,\beta}-\hat u_{i,\beta}\|_{\infty}\leq
\|u_{i,\beta}-g\|_{\infty}+\|\hat u_{i,\beta}-g\|_{\infty}\leq
C\beta^{-\frac{1}{2}}. $$
Concerning condition (c), it is clearly
satisfied in $A_i$, then it remains to check, for $\beta$ large
enough,  the inequality $\phi^z(x,z)+\beta(z-g)^2> 0$ holds
true outside $A_i$. For
$x\in\left(\Om_i\setminus(\Gamma)_{\frac{D}{2}}\right)\cap
\Gamma_{{\delta}}$ such an estimate can be performed using
estimates \eqref{gradu2}, \eqref{hessian2}, \eqref{ubetamod2} and
arguing as in the proof of Theorem \ref{calibration1}. Now let $(x,z)$ belong
to
$\left[\left(\Om_i\setminus(\Gamma)_{{D}}\right)\times\erre\right]\setminus
A_i$ and suppose
 also that $\hat u_{2,\beta}(x)+h_{\beta}(x)\leq z\leq\hat
u_{1,\beta}(x)-h_{\beta}(x)$ (the other cases would be analogous);
since $\phi^z(x,z)=\phi(x,\hat
 u_{i,\beta}+(-1)^{i}h_{\beta})\cdot
(-\nabla \hat
 u_{i,\beta}+(-1)^{i+1}\nabla h_{\beta},1)$ and observing that
$\phi(x,\hat u_{i,\beta}+(-1)^{i}h_{\beta})$ reduces to
$$
\left(2\nabla u_{\beta},\,|\nabla u_{\beta}|^2-\beta(z-g)^2+\beta
(u_{\beta}-z)^2\right),$$
we obtain
\begin{eqnarray*}
\phi^z(x,z)+\beta(z-g)^2&\geq& -|\nabla u_{\beta}||\nabla \hat u_{i,\beta}|
-2|\nabla u_{\beta}| |\nabla h_{\beta}|+\beta(u_{\beta}-z)^2\\
&\geq&
-|\nabla u_{\beta}||\nabla \hat u_{i,\beta}|-2|\nabla u_{\beta}| |\nabla
h_{\beta}|+\beta s_{\beta}^2;
\end{eqnarray*}
in the last expression the positive term $\beta
s_{\beta}^2$, which behaves like $\beta^{\frac{1}{2}}$ (see the definition of
$s_{\beta}$) dominates the negative ones, indeed these are either bounded or
of the same order of $|\nabla u_{\beta}|$ which is less or equal to the
order of $\beta^{\frac{1}{4}}$, thanks to \eqref{gradu2}: therefore for $\beta$
large  enough we get the desired inequality.

About condition (f) we first observe that   if $t_1$,
$t_2\in\erre$ and
 $x\in\left((\Gamma)_{D}\setminus(\Gamma)_{\frac{D}{2}}\right)\cap\Om$
then we obtain
$$
\left|\int_{t_1}^{t_2}\phi^x(x,z)\,
dz\right|\leq\qquad\qquad\qquad\qquad\phantom{llllllllllllllllllllll
lllllllllllllllllllllllllllllllllllllllllllllllllllllllllll}
$$
\vspace{-0.5cm}
\begin{eqnarray}
&\leq& \sum_{i=1}^2\left[\int_{ \hat u_{i,\beta}-h_{\beta}}^{
\hat u_{i,\beta}+h_{\beta}}\left(\hbox{$2\|\nabla
u_{i,\beta}\|_{\infty}+\frac{16}{\hat h_{\beta}} \|\nabla
u_{i,\beta}\|_{\infty}\left((-1)^i(z- u_{i,\beta})-\frac{\hat
h_{\beta}}{2}\right)^+ +4{|
u_{i,\beta}-z|}\|\nabla  v_i\|_{\infty} $} \right )\, dz\right]\nonumber\\
&\leq& \sum_{i=1}^2\left[4\|\nabla
u_{i,\beta}\|_{\infty}\varepsilon+ 4\varepsilon^2\|\nabla
v_i\|_{\infty}+\hbox{$\frac{16}{\hat h_{\beta}} \|\nabla
u_{i,\beta}\|_{\infty} \left[(u_{i,\beta}-\hat
u_{i,\beta})^2+\frac{\hat
h_{\beta}^2}{4}\right]$}\right]\nonumber\\
&\leq& \sum_{i=1}^2\left[4\|\nabla
u_{i,\beta}\|_{\infty}\varepsilon+ 4\varepsilon^2\|\nabla
v_i\|_{\infty}+\frac{16}{s_{\beta}} \|\nabla
u_{i,\beta}\|_{\infty}( u_{i,\beta}-\hat u_{i,\beta})^2+
8\|\nabla  u_{i,\beta}\|_{\infty}\varepsilon\right]\nonumber\\
&\leq& \sum_{i=1}^2\left[12\|\nabla
u_{i,\beta}\|_{\infty}\varepsilon+4\varepsilon^2\|\nabla
v_i\|_{\infty}+ C\beta^{-\frac{3}{4}}\right]\nonumber\\
&&\quad\left[\hbox{ the fact that $\frac{16}{s_{\beta}} \|\nabla
u_{i,\beta}\|_{\infty}( u_{i,\beta}-\hat
u_{i,\beta})^2\leq C\beta^{-\frac{3}{4}}$ follows from estimates}\right.\nonumber\\
&&\quad\left.\hbox{$\phantom\int$ \eqref{gradu2}, \eqref{vicine2},
\eqref{ubetamod2}, and the definition of
$s_{\beta}$}\right]\nonumber\\
&\leq&\frac{1}{2}\nonumber,
\end{eqnarray} if $\beta$ is large enough,
thanks to \eqref{epsilon2}.

If $x\in\Om_i\setminus(\Gamma)_{D}$ then we can estimate
$$
\left|\int_{t_1}^{t_2}\phi^x(x,z)\, dz\right|\leq 2
s_{\beta}\|\nabla  u_{1,\beta}\|_{\infty}+2
s_{\beta}\|\nabla  u_{2,\beta}\|_{\infty}\leq \frac{1}{2},
$$ by \eqref{gradu2} and the definition of $s_{\beta}$. Also
condition (f) is proved; since, by construction, $\phi$ has
vanishing normal component along $\partial\Om\times\erre$, the
theorem is completely proved.\qed

Now we can state a theorem which is the analogous of Theorem
\ref{cor1}.
\begin{theorem}\label{cor2}
Let $\Om$ as in Proposition \ref{BOOO} and
$\Gamma=\gamma_1\cup\dots\cup\gamma_k$  where for every
$j=1,\dots,k$   $\gamma_j$ is either a simple, connected, and
closed curve of class $C^{2,\alpha}$ contained in $\Om$ or a
connected curve with the same regularity outside a neighbourhood
of its endpoints (where it is supposed to be of class $C^3$),
which meets orthogonally $\partial\Om$ in two regular points;
suppose in addition that $\gamma_i\cap\gamma_j=\emptyset$ if
$i\neq j$. Then for every $g\in W^{1,\infty}(\Om\setminus\Gamma)$
discontinuous along $\Gamma$ and
 such that $g^+(x)-g^-(x)>S>0$ for every
$x\in\Gamma$, there exists $\beta_0>0$ depending on $\Gamma$, $S$,
and $\|g\|_{W^{1,\infty}}$, such that for $\beta\geq\beta_0$ the
solution $u_{\beta}$ of \eqref{ubeta'} is discontinuous along
$\Gamma$ ($S_{u_{\beta}}=\Gamma$) and it is the unique minimizer
of $F_{\beta,g}$ over $SBV(\Om)$.
\end{theorem}

\section{Gradient flow for  the Mumford-Shah functional}

In this section we are going to apply the previous results to the
study of the gradient flow of the Mumford-Shah functional
by the method of minimizing movements (see Section 2) with an initial datum
$u_0$ which
is regular outside a regular discontinuity set $\Gamma$: we
will show that, for an initial
interval of time,
 the discontinuity set does not move  while the  function evolves
 according
 to the heat equation.
Our main result is stated in the following theorem:
\begin{theorem}\label{movmin}
Let $\Om$ and $\Gamma$   be either as in Theorem \ref{cor1} or as
in Theorem \ref{cor2}. Suppose that $u_0$ is a function belonging
to $W^{2,\infty}(\Om\setminus\Gamma)$, discontinuous along
$\Gamma$, and  such that $u_0^+(x)-u_0^-(x)>S>0$ for every
$x\in\Gamma$ and $\partial_{\nu} u_0=0$ on
$\partial\Om\cup\Gamma$. Then there exists $T>0$ such that the
minimizing movement for the Mumford-Shah functional is unique in
$[0,T]$ and it is given by the function $u(x,t)$ satisfying
$$ S_{u(\cdot,t)}=\Gamma\qquad \forall t\in[0,T],$$
and
$$
\begin{cases}
\partial_t u=\Delta u &\text{in  $(\Om\setminus\Gamma)\times[0,T]$},\cr
\partial_{\nu}u=0 & \text{on  $\partial(\Om\setminus\Gamma)\times[0,T]$},\cr
 u(x,0)= u_0(x) & \text{in  $\Om\setminus\Gamma$}.
\end{cases}
$$
 \end{theorem}
{\sc Proof.}
For fixed $\delta>0$, let $v_{\delta}(t)$ be  the affine interpolation
 of the discrete function
\begin{eqnarray*}
v_{\delta}:\delta\enne &\to& H^1(\Om\setminus\Gamma) \\
v_{\delta }(\delta i) &\mapsto& v_{\delta,i},
\end{eqnarray*}
where $v_{\delta,i}$ is
 inductively defined as follows:
 \begin{equation}\label{vdeltai}
 \begin{cases}
 v_{\delta,0}= u_0,\cr
\cr
 v_{\delta,i} \hbox{  is the unique solution of}\cr
 \qquad\displaystyle{\min_{z\in
H^1(\Om\setminus\Gamma)}\int_{\Om\setminus\Gamma}|\nabla z|^2\
dx+\frac{1}{\delta}\int_{\Om\setminus\Gamma}|z-v_{\delta,i-1}|^2\ dx.}   \end{cases}
 \end{equation}

 {\sc Claim 1.} For every $T>0$, we have that
 $$v_{\delta}\to v\qquad \hbox{in }
 L^{\infty}([0,T];L^{\infty}(\Om\setminus\Gamma))\ \hbox{as } \delta\to 0,$$
 where $v$ is the solution  of
\begin{equation}\label{vu}
\begin{cases}
\partial_t v=\Delta v &\text{in  $(\Om\setminus\Gamma)\times[0,T]$},\cr
\partial_{\nu}v=0 & \text{on  $\partial(\Om\setminus\Gamma)\times[0,T]$},\cr
 v(x,0)= u_0(x) & \text{in  $\Om\setminus\Gamma$}.
\end{cases}
\end{equation}

We will show that the functions $(v_{\delta})_{\delta>0}$ are equibounded
in $C^{0,1}([0,T]; L^{\infty}(\Om\setminus\Gamma))$: since
it is well known that, for every $T>0$, $v_{\delta}\to v$ in
   $L^{\infty}([0,T];L^{2}(\Om\setminus\Gamma))\ \hbox{as } \delta\to 0$
 (see for example \cite{Amb0}),
the a priori estimate in the $C^{0,1}$-norm (via Ascoli-Arzel\`a Theorem)
 will give the thesis of Claim 1.
 First of all we will show that
\begin{equation}\label{laplaciani}
\|\Delta v_{i,\delta}\|_{\infty}\leq \|\Delta u_0\|_{\infty}\qquad
\forall\delta>0,\, \forall i\in\enne.
\end{equation}
We first prove it for $v_{\delta,1}$: if
$\varepsilon\geq\|\Delta u_0\|_{\infty}/\beta$, then
$v_1:=u_0+\varepsilon$ and $v_2:=u_0-\varepsilon$ satisfy:
$$
\begin{cases}
\Delta v_1\leq\beta(v_1-u_0)&\text{in $\Om\setminus\Gamma$}\\
\partial_{\nu}v_1=0 &\text{on $\partial(\Om\setminus\Gamma)$,}
\end{cases}
\begin{cases}
\Delta v_2\geq\beta(v_2-u_0)&\text{in $\Om\setminus\Gamma$}\\
\partial_{\nu}v_2=0 &\text{on $\partial(\Om\setminus\Gamma)$,}
\end{cases}
$$
that is $v_1$ and $v_2$ are a supersolution and a subsolution respectively
of the problem solved by $v_{1,\delta}$. This implies that
$$
\|v_{1,\delta}-u_0\|_{\infty}\leq\frac{\|\Delta
u_0\|_{\infty}}{\beta}
$$
which is equivalent to
$$
\|\Delta v_{1,\delta}\|_{\infty}\leq \|\Delta u_0\|_{\infty}.
$$
By the same argument we can prove that
$$
\|\Delta v_{i,\delta}\|_{\infty}\leq\|\Delta
v_{i-1,\delta}\|_{\infty}\qquad\forall i\geq i-1 $$
and so \eqref{laplaciani} follows by induction on $i$.

By a standard truncation argument, one can prove also
 that
\begin{equation}\label{elleinfinito}
\|v_{\delta,i}\|_{\infty}\leq
\|u_0\|_{\infty}\qquad\forall\delta>0,\ \forall i\in\enne.
\end{equation}
Then for $s$, $t>0$, using Claim 1, we can estimate
$$
\|v_{\delta}(t)-v_{\delta}(s)\|_{\infty}\leq
\int_s^t\|(v_{\delta})'(\xi)\|_{\infty}\ d\xi\leq
\int_s^t\sup_i\|\Delta v_{\delta,i}\|_{\infty}\ d\xi\leq\|\Delta
u_0\|_{\infty}|t-s|;$$ this, together with (\ref{elleinfinito})
concludes the proof of Claim 1.

As a consequence of \eqref{laplaciani}, by the well-known Calderon-Zygmund
estimates, we get the existence of a constant $C$ such that
\begin{equation}\label{gradienti}
\|\nabla v_{i,\delta}\|_{\infty}\leq C\|\Delta v_{i,\delta}\|_{\infty}\leq
C\|\Delta u_0\|_{\infty}\qquad \forall\delta>0,\, \forall
i\in\enne.
\end{equation}
 It is well known (see, for example, \cite{LU})
that
$$
v(t)\to u_0\qquad \hbox{in }L^{\infty}(\Om\setminus\Gamma)\hbox{
as } t\to 0^+;
$$
 therefore, by our assumption on $u_0$,
 for every $0<c<S$ we can find $T_c>0$ such that
\begin{equation}\label{1401}
\inf_{x\in\Gamma}|v^+(x,t)-v^-(x,t)|>{c}\qquad \forall
t\in[0,T_c],
\end{equation}
and therefore, by Claim 1, we can choose $\delta_0>0$ such that
 \begin{equation}\label{1402}
\inf_{x\in\Gamma}|v^+_{\delta}(t,x)-v^-_{\delta}(t,x)|>\frac{c}{2}\qquad
\forall t\in[0,T_c],\   \forall\delta\leq\delta_0.
\end{equation}
We recall now that, by Theorems \ref{cor1} and \ref{cor2}, there
exists   $\overline{\beta}$ such that, for every function $g\in
W^{2,\infty}(\Om\setminus\Gamma)$ satisfying
 \begin{equation}\label{1403}
\|\nabla g\|_{\infty}\leq C\|\Delta u_0\|_{\infty}\qquad
\inf_{x\in\Gamma}|g^+(x)-g^-(x)|>\frac{c}{2},
\end{equation}
 where $C$ is the constant appearing in \eqref{gradienti},
and for every $\beta\geq\overline{\beta}$, the function
$u_{\beta,g}$ solution of \eqref{ubeta'}, minimizes the functional
$F_{\beta,g}$ over $SBV(\Om)$.

{\sc Claim 2.} For every
$\delta\leq\delta_0\land(\overline{\beta})^{-1}$
 the $\delta$-approximate evolution $u_{\delta}(t)$ (see the end of
 Section 2 for the definition) coincides in the interval $[0,T_c]$
  with the function $v_{\delta}(t)$.

Clearly it is enough to show  that
$$
v_{\delta,i}=u_{\delta,i}\qquad \hbox{for }i=0,\dots, \left[\frac{T_c}{\delta}\right],
$$
and this can be done by induction on $i$:
indeed for $i=0$ the identity is trivial, and suppose it true for $i-1$
 (for $i\leq\left[\frac{T_c}{\delta}\right]$ ); this means in  particular
 (by \eqref{gradienti} and by (\ref{1402})) that $g=u_{\delta,i-1}$ satisfies
(\ref{1403})
  and so, being $\frac{1}{\delta}>\overline{\beta}$,  we have
$$ u_{\delta,i}=u_{\frac{1}{\delta},u_{\delta,i-1}}=v_{\delta,i}.$$
 Claim 2 is proved and the thesis of the theorem is now evident.
\qed

\section*{Acknowledgements}
I am grateful to Gianni Dal Maso for  drawing my
attention on this problem and for interesting
discussions.


\begin{thebibliography}{99}

\bibitem{Alb-Bou-DM} Alberti G., Bouchitt\'e G., Dal Maso G.:
The calibration method for the Mumford-Shah functional.  {\it
C. R. Acad. Sci. Paris S\'er. I Math.} {\bf 329} (1999), 249-254.

\bibitem{Alb-Bou-DM2} Alberti G., Bouchitt\'e G., Dal Maso G.:
The calibration method for the Mumford-Shah functional and free
discontinuity problems. Preprint  SISSA, Trieste, 2001.

\bibitem{Amb0} Ambrosio L.: Movimenti minimizzanti. {\it Rend. Accad. Naz.
Sci. XL Mem. Mat. Sci. Fis. Natur.} {\bf 113} (1995), 191-246.

\bibitem{Amb} Ambrosio L.: A compactness theorem for a new class of
variational problems. {\it Boll. Un. Mat. It.\/} {\bf 3-B} (1989),
857-881.

\bibitem{Amb-Fus-Pal} Ambrosio L., Fusco N., Pallara D.:
Special Functions of Bounded Variation and Free-Discontinuity
Problems. Oxford University Press, Oxford, 2000.

\bibitem{Bonnet} Bonnet A.: On the regularity of edges in image
segmentation. {\it Ann. Inst. H. Poincar\'e, Anal. non
lin\'eaire.\/} {\bf 13} (1996), 485-528.

\bibitem{Cha-Dov} Chambolle A., Doveri F.: Minimizing movements of
the Mumford-Shah energy. {\it Discrete  Contin. Dynam. Systems\/}
{\bf 3} (1997), 153-174.

\bibitem{DM-M-M} Dal Maso G., Mora M. G., Morini M.:
Local calibrations for minimizers of the Mumford-Shah functional
with rectilinear discontinuity set.  {\it J. Math. Pures Appl.\/}
{\bf 79} (2000), 141-162.

\bibitem{DG-Amb} De Giorgi E., Ambrosio L.: Un nuovo funzionale
del calcolo delle variazioni. {\it Atti Accad. Naz. Lincei Rend.
Cl. Sci. Fis. Mat. Natur. (8)} {\bf 82} (1988), 199-210.


\bibitem{DZ} Delfour M. C., Zol\'esio J. P.: Shape Analysis via
distance functions. {\it J. Functional Analysis \/} {\bf 123}
(1994), 129-201.

\bibitem{Gob} Gobbino M.: Gradient flow for the one-dimensional
Mumford-Shah functional. {\it Ann. Scuola Norm. Sup. Pisa Cl. Sci.
(4) \/} {\bf 27} (1998), 145-193.

\bibitem{GRI1} Grisvard P.: Majorations en norme du maximum de la
r\'esolvante du laplacien dans un polygone. Nonlinear partial
differential equations and their applications. Coll\`ege de France
Seminar, Vol. XII (Paris, 1991-1993), 87-96, Pitman Res.

\bibitem{GRI2} Grisvard P.: Elliptic Problems in Nonsmooth
Domains. Monographs and Studies in Mathematics, 24. Pitman
(Advanced Publishing Program), Boston, MA, 1985.

\bibitem{LU} Lunardi A.: Analytic Semigroups and Optimal
Regularity in Parabolic Problems. Progress in Nonlinear
Differential Equations and Their Applications {\bf 16},
Birkh\"auser Verlag, Basel, 1995.

\bibitem{Circe} Mora M. G.: Local calibrations for minimizers of
the Mumford-Shah functional with a triple junction. Preprint SISSA, Trieste,
2001.

\bibitem{mora-morini-99} Mora M. G., Morini M.: Functional
depending on curvatures with constraints. {\it Rend. Sem. Mat.
Univ. Padova \/} {\bf 104} (2000), 173-199.

\bibitem{M-M} Mora M. G., Morini M.: Local calibrations for
minimizers of the Mumford-Shah functional with a regular
discontinuity set. To appear on {\it Ann. Inst. H. Poincar\'e,
Anal. non lin\'eaire.\/}

\bibitem{Mum-Sha1} Mumford D., Shah J.:
Boundary detection by minimizing functionals, I. {\it Proc. IEEE
Conf. on Computer Vision and Pattern Recognition (San Francisco,
1985)\/}.

\bibitem{Mum-Shah} Mumford D., Shah J.:
Optimal approximation by piecewise smooth functions and associated
variational problems. {\it Comm. Pure Appl. Math.\/} {\bf 42}
(1989), 577-685.

\bibitem{Rich} Richardson T. J.: Limit theorems for a variational
problem arising in computer vision. {\it Ann. Scuola Norm. Sup.
Pisa Cl. Sci. (4) \/} {\bf 19} (1992), 1-49.

\end{thebibliography}
\end{document}